\documentclass[12pt]{article}
\usepackage{latexsym,amsmath,amssymb}
\usepackage{graphicx}
\begin{document}

\title{{\large A NEW UPPER BOUND ON THE ACYCLIC CHROMATIC INDICES OF PLANAR
GRAPHS}}

\author{Weifan Wang\thanks{Research supported partially by
NSFC(No.11071223) and ZJNSF(No.Z6090150); Corresponding author.
Email: wwf@zjnu.cn.}, \ Qiaojun Shu\\
\bigskip
 \normalsize Department of Mathematics,  Zhejiang Normal University, Jinhua 321004, China\\
 Yiqiao Wang\\
 \normalsize Academy of Mathematics and Systems Science, Beijing 100080,
 China}

\date{ }

\maketitle
\newtheorem{question}{Question}
\newtheorem{theorem}{Theorem}
\newtheorem{lemma}{Lemma}
\newtheorem{problem}{Problem}
\newtheorem{claim}{Claim}
\newtheorem{corollary}[theorem]{Corollary}
\newtheorem{definition}{Definition}
\newtheorem{proposition}[theorem]{Proposition}
\newtheorem{conjecture}{Conjecture}
\newcommand{\qed}{\hfill $\Box$ }
\newcommand{\proof}{\noindent{\bf Proof.}\ \ }
\baselineskip=15pt

\parindent=0.7cm

%%%%%%%%%%%%%%%%%%%%%%%%%%%%%%%%%%%%%%%%%%%%%%%%%%%%%%%%%%%%%%%%%
\begin{abstract}

\baselineskip=16pt

An acyclic edge coloring of a graph $G$ is a proper edge coloring
such that no bichromatic cycles are produced. The acyclic chromatic
index $a'(G)$ of $G$ is the smallest integer $k$ such that $G$ has
an acyclic edge coloring using $k$ colors.  It was conjectured that
  $a'(G)\le \Delta+2$ for any simple graph $G$ with maximum
degree $\Delta$. In this paper, we prove that if $G$ is a planar
graph, then $a'(G)\leq\Delta +7$. This improves a result by
 Basavaraju et al. [{\em Acyclic edge-coloring of planar graphs},   SIAM
J. Discrete Math.,  25 (2011), pp. 463-478], which says that every
planar graph $G$ satisfies  $a'(G)\leq\Delta +12$.


\medskip


\end{abstract}

%%%%%%%%%%%%%%%%%%%%%%%%%%%%%%%%%%%%%%%%%%%%%%%%%%%%%%%%%%%%%%%%
\section{Introduction}
%%%%%%%%%%%%%%%%%%%%%%%%%%%%%%%%%%%%%%%%%%%%%%%%%%%%%%%%%%%%%%%%

\baselineskip=18pt

Only simple graphs are considered in this paper. Let $G$ be a graph
with vertex set $V(G)$ and edge set $E(G)$. A {\em proper edge
$k$-coloring} is a mapping $c:E(G) \to \{1,2,\ldots ,k\}$ such that
any two adjacent edges receive different colors.  The {\em chromatic
index} $\chi'(G)$ of $G$ is the smallest integer $k$ such that $G$
is edge $k$-colorable. A proper edge $k$-coloring $c$ of $G$ is
called {\em acyclic} if there are no bichromatic cycles in $G$,
i.e., the union of any two color classes induces a subgraph of $G$
that is a forest. The {\em acyclic chromatic index} of $G$, denoted
by $a'(G)$, is the smallest integer $k$ such that $G$ is acyclically
edge $k$-colorable.

Let $\Delta(G)$ ($\Delta$ for short) denote the maximum degree of a
graph $G$ and $g(G)$ denote the girth of $G$, i.e., the length of a
shortest cycle in $G$. By Vizing's Theorem \cite{Vizing64}, $\Delta
\le \chi'(G)\le \Delta +1$. Thus, it is obvious that $a'(G)\ge
\chi'(G)\ge \Delta$. Fiam${\rm \check{c}}$ik \cite{F78},  and later
Alon, Sudakov and Zaks \cite{ABZ01},  made independently the
following conjecture:
\begin{conjecture}\label{con1}
For any graph $G$, $a'(G)\le \Delta+2$.
\end{conjecture}

Using probabilistic method, Alon, McDiarmid and Reed \cite{AMR91}
proved that $a'(G)\le 64\Delta$ for any graph $G$. Molloy and Reed
\cite{MR98} improved this bound to that  $a'(G)\le 16\Delta$. Alon,
Sudakov and Zaks \cite{ABZ01} proved that there is a constant $c$
such that $a'(G)\le \Delta+2$ for a graph $G$ whenever $g(G)\ge
c\Delta\log \Delta$. They also confirmed Conjecture \ref{con1} for
almost all $\Delta$-regular graphs. N$\check{\rm e}$set$\check{\rm
r}$il and Wormald \cite{NNCW05} gave an improvement to this result
by showing that $a'(G)\le \Delta+1$ for a random $\Delta$-regular
graph $G$. The acyclic edge coloring of some special classes of
graphs was also considered, including subcubic graphs \cite{BLSC08,
Sku04}, graphs with maximum degree $4$ \cite{MLSC09}, outerplanar
graphs \cite{HWLL10, MNCRS07}, series-parallel graphs
\cite{HWLL09,WS}, and planar graphs \cite{BF10, DX10, HLW11, SW10,
WSWW11, SWW11, YHLLX09}.

%In particular,  Conjecture \ref{con1} is confirmed for some special
%classes of planar graphs, see \cite{HLW11, SW11, WSW10, SWW11}.

Suppose that $G$ is a planar graph. Fiedorowicz et al.\,\cite{FHN08}
proved that $a'(G)\leq 2\Delta+ 29$ and Hou et al.\,\cite{HWLL09}
proved  that $a'(G)\leq \max\{2\Delta- 2, \Delta + 22\}$. The
currently best known upper bound on the acyclic chromatic index of a
planar graph $G$ is that $a'(G)\le \Delta + 12$  by Basavaraju et
al.\,\cite{BCCHM11}. In this paper, we will improve this result by
replacing $12$ by $7$.

Before showing the main result, we need to introduce some
definitions and notations.

Given a graph $G$, let $d_G(v)$ (or simply $d(v)$) denote the degree
of a vertex $v$ in $G$. A vertex of degree $k$ (at least $k$, at
most $k$) is called a $k$-{\em vertex} ({$k^+$-vertex,
$k^-$-vertex}). For $k\ge 1$, let $n_k(v)$ ($n_{k^{+}}(v)$,
$n_{k^{-}}(v)$) denote the number of $k$-vertices ($k^{+}$-vertices,
$k^{-}$-vertices) adjacent to a vertex $v$ in $G$. Let $\delta(G)$
denote the minimum degree of $G$.

A {\it plane} graph is a particular drawing in the Euclidean plane
of a certain planar graph. For a plane graph $H$, we denote its face
set by $F(H)$ and the degree of a face $f\in F(H)$ is denoted by
$d_H(f)$ (or simply $d(f)$). Similarly, we can define a $k$-{\em
face}, a $k^+$-{\em face}, and a $k^-$-{\em face}. Furthermore, for
$f\in F(H)$, we use $b(f)$ to denote the boundary walk of $f$ and
write $f= [u_1u_2\ldots u_n]$ if $u_1, u_2, \ldots, u_n$ are the
vertices on $b(f)$ enumerated in some order, where $u_i$ may be
identical to $u_j$ for some $i, j \in \{1, 2, \ldots, n\}$ when $G$
contains cut vertices.

\medskip

%%%%%%%%%%%%%%%%%%%%%%%%%%%%%%%%%%%%%%%%%%%%%%%%%%%%%%%%%%%%%%
\section{A structural lemma}
%%%%%%%%%%%%%%%%%%%%%%%%%%%%%%%%%%%%%%%%%%%%%%%%%%%%%%%%%%%%%%

\medskip

In this section, we establish a structural lemma, which plays an
important role in the proof of the main result in Section 3.

\begin{lemma}\label{planar}

Let $G$ be a $2$-connected planar graph with $\Delta\ge 5$. Then $G$
contains one of the following configurations {\bf
(A$_{1}$)-(A$_{4}$)}, as shown in Fig.\,1:

{\bf (A$_{1}$)}\ A path $uvw$ with $d(v)= 2$ and $d(u)\le 9$.

{\bf (A$_{2}$)}\ A vertex $u$ with $n_{2}(u)\geq 1$ and
$n_{8^{-}}(u)\geq d(u)-8$. Suppose that $u_{1}, u_{2}, \ldots, $
$u_{d(u)-1}$, $v$ are the neighbors of $u$ such that $d(u_{1})\geq
d(u_{2})\geq \cdots \geq d(u_{d(u)-1})\geq d(v)= 2$. Let $w$ be the
neighbor of $v$ different from $u$. For $1\le i\le d(u)-1$, if
$d(u_i)=2$, then we use $x_i$ to denote the neighbor of $u_i$
different from $u$. Then at least one of the following cases holds:

\quad  {\bf (A$_{2.1}$)}\ $n_{8^{-}}(u)\geq d(u) - 7$;

\quad  {\bf (A$_{2.2}$)}\ $n_{8^{-}}(u) = d(u) - 8$, and
$n_{2}(u)\geq d(u)-9$.

  {\bf (A$_{3}$)}\ A $3$-vertex $u$ is adjacent to a vertex
$v$ such that one of the following holds:

\quad  {\bf (A$_{3.1}$)}\ $d(v)\leq 8$;

\quad  {\bf (A$_{3.2}$)}\ $d(v)= 9$ and $uu_2, vu_2 \in E(G)$;

\quad  {\bf (A$_{3.3}$)}\ $d(v)= 10$, $n_{5^-}(v)\geq 5$ and $uu_i,
vu_i \in E(G)$, $i= 1, 2$.

 {\bf (A$_{4}$)}\ Suppose that a vertex $v$ is adjacent to
$u, v_2, \cdots, v_{d(v)}$ such that $d(u)\leq d(v_2) \leq \cdots
\leq d(v_{d(v)})$. Then at least one of the following cases holds:

\quad {\bf ($A_{4.1}$)}\ $d(v)= 4$, $4 \leq d(u) \leq 7$, and $d(u)
+ d(v_2) \leq 17$;

\quad {\bf ($A_{4.2}$)}\ $d(v)= 5$, $4 \leq d(u) \leq 6$, and $d(u)
+ d(v_2)+ d(v_3) \leq 18$, or $d(u)= d(v_2)= 6$, $d(v_i)= 7$, $i= 3,
4, 5$ and $uv_5\in E(G)$.

\end{lemma}

\begin{figure}[hh]
\begin{center}

\includegraphics[width=4.5 in]{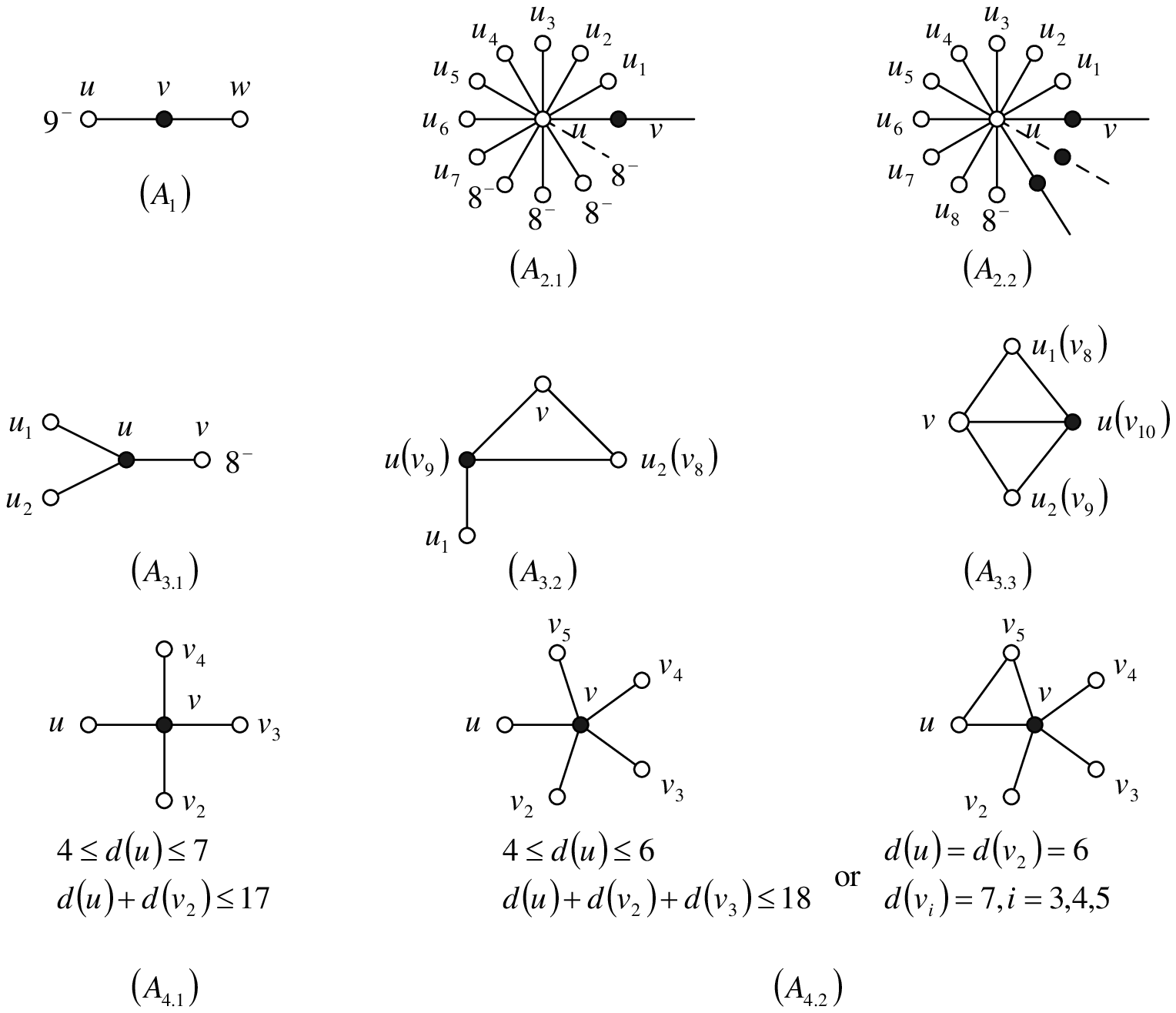}\\

\bigskip

{\rm Fig.\,1: \ Configurations ($A_{1}$)-($A_{4}$) in Lemma
\ref{planar}.}
\end{center}
\end{figure}

\noindent{\bf Remarks.}\  In Fig.\,1, vertices marked solid points
 have no edges of $G$ incident to them other than those shown,
 whereas the other vertices have edges connected to other vertices of
 $G$ not in the configuration.\\

\proof\ Assume to the contrary that $G$ contains none of the
configurations ($A_{1}$)-($A_{4}$). Since $G$ is $2$-connected, it
follows that $\delta(G)\ge 2$. Let $G'$ be the graph obtained by
removing all the $2$-vertices of $G$ and $H$ is a component of $G'$.
Then $H$ is a connected planar graph and for any vertex $v\in V(H)$,
$v$ is of degree at least $3$ in $G$.

In what follows, we assume that $H$ is embedded in the plane. Let
$v\in V(H)$. Then $v\in V(G)$.  The degree of $v$ in $G$ is simply
denoted by $d(v)$, whereas the degree of $v$ in $H$ is denoted by
$d_H(v)$. For $k\ge 1$, let $n'_k(v)$ ($n'_{k^{+}}(v)$,
$n'_{k^{-}}(v)$) denote the number of $k$-vertices
($k^{+}$-vertices, $k^{-}$-vertices) adjacent to  $v$ in $H$, while
$n_k(v)$ ($n_{k^{+}}(v)$, $n_{k^{-}}(v)$) denote the number of
$k$-vertices ($k^{+}$-vertices, $k^{-}$-vertices) adjacent to  $v$
in $G$. The number of $k$-vertices ($k^{+}$-vertices) lying on the
boundary of a face $f\in F(H)$ is  denoted by $n_k(f)$
($n_{k^{+}}(f)$). Similarly, let $m_{k}(v)$ denote the number of
$k$-faces incident to  $v$ in $H$, and let $\delta(f)$ denote the
minimum degree of vertices on $b(f)$.

Since $G$ contains no ($A_{1}$), there is no $9^-$-vertex adjacent
to a $2$-vertex, i.e., if $u$ is a vertex with $3\leq d(u)\leq 9$,
then $d_H(u)= d(u)$. Similarly, if $d(u)\geq 10$ and $n_2(u)= 0$,
then $d_H(u)= d(u)$. Since $G$ contains neither ($A_{2.1}$) nor
($A_{2.2}$), every $10^+$-vertex $u$ is adjacent to at most
($d(u)-9$) $2$-vertices, which implies that if $d(u)\geq 10$ and
$n_2(u)\geq 1$, then $d_H(u)\geq 9$, $d_H(u)= d(u)- n_{2}(u)$.

\begin{claim}\label{minimum}
Let $u\in V(H)$.

  {\bf (a)}\ $\delta(H)\geq 3$.

  {\bf (b)}\  If $d_H(u)\leq 8$, then $d(u)= d_H(u)$.

  {\bf (c)}\  $n'_{i}(u)= n_{i}(u)$ and $n_{i^-}(u)= n_2(u) +
n'_{i^-}(u) = d(u)- d_H(u) + n'_{i^-}(u)$ for any $3 \leq i \leq 8$.
In particular, $n'_{8^-}(u) = \sum\limits_{3\leq i\leq 8}n_i'(u)=
\sum\limits_{3\leq i\leq 8}n_i(u)= n_{8^-}(u)- n_2(u)$.
\end{claim}
\proof  It is easy to see that (c) holds trivially. To prove (a),
assume that $H$ contains a $2^-$-vertex $u$. Then $d(u)\geq 10$,
$n_{8^{-}}(u)\geq n_{2}(u) = d(u)- d_H(u)\geq d(u) - 2$ and hence
($A_{2.1}$) is contained in $G$, a contradiction.

To prove (b), assume that $d_H(u)\leq 8$ and $d(u)\neq d_H(u)$. It
is easy to see that $d(u)\geq 10$, and $n_{2}(u) = d(u)- d_H(u)\geq
d(u)- d_H(u)\geq d(u)- 8$, and hence ($A_{2.1}$) or ($A_{2.2}$) is
contained in $G$, also a contradiction. \qed

\begin{claim}\label{9-vertex}
If $d_H(u)= 9$ and $n'_{8^-}(u)\geq 1$, then $d(u)= 9$. In other
words, if $n'_{8^-}(u)\geq 1$, then $d_H(u)= d(u)\leq 9$, or
$d_H(u)\geq 10$.
\end{claim}
\proof If $d(u)\neq 9$, then $n_{2}(u)= d(u)- 9$, $n_{8^{-}}(u) \geq
n_{2}(u) + n'_{8^{-}}(u)\geq n_{2}(u) + 1 = d(u)- d_H(u)+ 1= d(u) -
8$ and hence ($A_{2.2}$) is contained in $G$. \qed

\begin{claim}\label{10+345}
If $d_H(u)\geq 10$ and $d(u)\neq d_H(u)$, then $n_3(u) + n_4(u) +
n_5(u)\leq d_H(u) - 8$.
\end{claim}
\proof Otherwise, assume that $n_3(u) + n_4(u) + n_5(u)\geq d_H(u) -
7$. Then, $n_{8^{-}}(u)\geq n_{2}(u) + n_{3}(u) + n_{4}(u) +
n_{5}(u)= d(u)- d_H(u) + n_{3}(u) + n_{4}(u) + n_{5}(u) \geq d(u)-
d_H(u) + d_H(u) - 7 = d(u)- 7$ and hence ($A_{2.1}$) is contained in
$G$. \qed

\begin{claim}\label{3-face}
If $f = [uvw]$ with $d_H(v)= 3$, then $d_H(u)\geq 10$, $d_H(w)\geq
10$ and $n_{10^+}(f)= 2$.
\end{claim}
\proof Since $G$ contains neither ($A_{3.1}$) nor ($A_{3.2}$), we
have $d(u)\geq 10$ and $d(w)\geq 10$. Obviously, $n'_{8^-}(u)\geq
n'_{3}(u)\geq 1$ and $n'_{8^-}(w)\geq n'_{3}(w)\geq 1$. Then, by
Claim \ref{9-vertex}, $d_H(u)\geq 10$, $d_H(w)\geq 10$ and
$n_{10^+}(f)= 2$. \qed

\begin{claim}\label{10+10}
If $3 \leq d(u)\leq 8$, then $n'_{10^+}(u)= n_{10^+}(u)$.
\end{claim}
\proof  Note that $d_H(w)\leq d(w)$ for any $w\in V(H)$. Let $v$ be
a vertex of degree at least $10$ in $G$ and $uv\in E(G)$. Since
$uv\in E(H)$, $n'_{8^-}(v)\geq n_{3}(v)\geq 1$ and $d_H(v)\geq 10$
by Claim \ref{9-vertex}. \qed

\begin{claim}\label{4-vertex}
Let $u$ be a $4$-vertex. If $n_{7^-}(u)\geq 1$, then $n'_{10^+}(u)=
3$. Otherwise, $n'_{8^+}(u)= 4$. \qed
\end{claim}
\proof Since $G$ contains neither ($A_{1}$) nor ($A_{3.1}$),
$n_{3^-}(u)= 0$. If $n_{7^-}(u)\geq 1$, then $n_{i}(u)\geq 1$ for
some $4\leq i\leq 7$. Since $G$ contains no ($A_{4.1}$), we have
$n_{11^+}(u)= 3$ and $n'_{10^+}(u)= 3$ by Claim \ref{10+10}.
Otherwise, $n_{8^+}(u)= 4$ and $n'_{8^+}(u)= 4$. \qed

\begin{claim}\label{5-vertex}
Let $u$ be a $5$-vertex with $n_{4}(u)\geq 1$. If $n_{5^-}(u)\geq
2$, then $n'_{10^+}(u)= 3$.
\end{claim}
\proof Since $G$ contains neither ($A_{1}$) nor ($A_{3.1}$),
$n_{3^-}(u)= 0$. Since $G$ contains no ($A_{4.2}$), we have
$n_{10^+}(u)= 3$ and $n'_{10^+}(u)= 3$ by Claim \ref{10+10}. \qed

To derive a contradiction, we make use of the discharging method.
First, by Euler's formula $|V(H)|-|E(H)|+|F(H)| = 2$ and the
relation $\sum\limits_{u\in V(H)}d_H(u) = \sum\limits_{f\in
F(H)}d(f) = 2|E(H)|$, we can derive the following identity.
\begin{equation}\label{eqa1}
\sum_{u\in V(H)}(2d_H(u)-6) + \sum_{f\in F(H)}(d(f)-6) = - 12.
\end{equation}

Next, we define an initial weight function $w$ by $w(u) = 2d_H(u)-6$
for $u\in V(H)$ and $w(f) = d(f) - 6$ for $f\in F(H)$. It follows
from (\ref{eqa1}) that the total sum of weights is equal to $-12$.
In what follows, we will define some discharging rules and
redistribute weights accordingly. Once the discharging is finished,
a new weight function $w'$ is produced. However, the total sum of
weights is kept fixed when the discharging is in process.
Nevertheless, we can show that $w'(x)\ge 0$ for all $x\in V(H)\cup
F(H)$. This leads to the following obvious contradiction
$$
0\le \sum\limits_{x\in V(H)\cup F(H)}w'(x) = \sum\limits_{x\in
V(H)\cup F(H)} w(x) = -12<0
$$
and hence demonstrates that no such counterexample can exist.

A $3$-face $f$ is called {\em bad} if $f$ is incident to a
$3$-vertex, i.e., $\delta(f)= 3$.

Let $y$ be a vertex in $H$ and $f= [\cdots xyz \cdots]$ be a face
incident to $x, y, z$ with $d_H(x)\geq d_H(z)$.  Let $\tau(
y\rightarrow f)$ denote the amount of weight transferred from $y$ to
$f$ according to the following defined  discharging rules:
\medskip

\noindent{\bf (R1)}\ $d_H(y)= 4$. If $n_{7^-}(y)= 0$, then
$\tau(y\to f)= \frac{1}{2}$. If $n_{7^-}(y)= 1$, then
     $$\tau(y\to f)=
\begin{cases} \frac{4}{5},\ \text{if}\ d_H(z)\leq 7;\\
              \frac{1}{5},\ \text{if}\ d_H(z)\geq 8.\\
\end{cases}$$

\noindent{\bf (R2)}\ $d_H(y)= 5$. If $n_4(y)\geq 1$, then $\tau(y\to
f)= \frac{4}{5}$. Otherwise, $n_4(y)= 0$, we carry out the following
subrules:

$\bullet$\   If $4\leq d(f)\leq 5$ and $d_H(z)\leq 8$, then
$\tau(y\to f)= \frac{1}{2}$£»

$\bullet$\ If $d(f)= 3$, or $4\leq d(f)\leq 5$ and $d_H(z)\geq 9$,
then

$$\tau(y\to f)=
\begin{cases} \frac{7}{5},\ \text{if}\ d_H(x)= d_H(z)= 5;\\
              \frac{6}{5},\ \text{if}\ d_H(x)=6\ \text{and}\ d_H(z)= 5;\\
              \frac{13}{14},\ \text{if}\ d_H(x)=7\ \text{and}\ d_H(z)= 5;\\
              1,\ \ \text{if}\    d_H(x)=8\ \text{and}\ d_H(z)= 5;\\
              \frac{11}{12},\ \text{if}\   d_H(x)\geq 9\ \text{and}\ d_H(z)= 5;\\

              1,\ \text{if}\ d_H(x)= d_H(z)= 6;\\
              \frac{6}{7},\ \text{if}\ d_H(x)= 7\ \text{and}\ d_H(z)= 6;\\
              \frac{3}{4},\ \text{if}\ d_H(x)= 8\ \text{and}\ d_H(z)= 6;\\
              \frac{2}{3},\ \text{if}\ d_H(x)\geq 9\ \text{and}\ d_H(z)= 6;\\

              \frac{5}{7},\ \text{if}\ d_H(x)= d_H(z)= 7;\\
              \frac{9}{14},\ \text{if}\ d_H(x)\geq 8\ \text{and}\ d_H(z)= 7;\\

              \frac{1}{2},\ \text{if}\ d_H(z)= 8;\\
              \frac{1}{3},\ \text{if}\ d_H(z)\geq 9.\\

\end{cases}$$

\noindent{\bf (R3)}\ $d_H(y)= k$ and $6\leq k \leq 9$. Then
$$\tau(y\to f)=
\frac{2k-6}{k}= \begin{cases} \frac{4}{3},\ \text{if}\ k= 9;\\
              \frac{5}{4},\ \text{if}\ k= 8;\\
              \frac{8}{7},\ \text{if}\ k= 7;\\
              1,\ \text{if}\ k= 6.\\
\end{cases}$$

\noindent{\bf (R4)}\ $d_H(y)\geq 10$. If $4\leq d(f)\leq 5$,  or
$d(f)= 3$ and $d_H(z)\geq 6$, then $\tau(y\to f)=1$. Otherwise, $f=
[xyz]$ with $d_H(z)\leq 5$. Then $$\tau(y\to f)=
\begin{cases} \frac{3}{2},\ \text{if}\ d_H(z)= 3;\\

             \frac{7}{5},\ \text{if}\ d_H(z)= 4\ \text{and}\ n_{7^-}(z)= 1;\\
             \frac{5}{4},\ \text{if}\ d_H(z)= 4\ \text{and}\ n_{7^-}(z)= 0;\\

              \frac{11}{10},\ \text{if}\ d_H(z)= 5,\ n_4(z)\geq 1,\ \text{and}\ d_H(x)\geq 10;\\
              \frac{5}{4},\ \text{if}\ \ d_H(z)= 5,\ n_4(z)\geq 1, \ \text{and}\ 6\leq d_H(x)\leq 9;\\
              \frac{7}{5},\ \text{if}\ d_H(x)= d_H(z)= 5\ \text{and}\ n_4(z)\geq 1\ \text{or} \ n_4(x)\geq 1;\\
              \frac{4}{3},\ \text{if}\ d_H(z)= 5,\ n_4(z)= 0\ \text{and}\ n_4(x)= 0\ \text{if}\ d_H(x)= 5.\\

\end{cases}$$

It remains to inspect that $w'(x)\ge 0$ for all $x\in V(H)\cup
F(H)$. Let $f\in F(H)$. We consider several cases as follows:

\medskip

\noindent{\bf Case 1.}\  $d(f) = 3$.

Let $f= [uvw]$. Then  $w(f) = d(f)- 6 = -3$. We have some subcases,
depending on the size of $\delta(f)$.

{\bf (1.1)}\  Assume that $\delta(f) = 3$. Then $n_{10^+}(f) = 2$ by
Claim \ref{3-face} and  (R4), each $10^+$-vertex in $b(f)$ gives
$\frac{3}{2}$ to $f$. Hence, $w'(f) = -3 + 2\times \frac{3}{2} = 0$.

{\bf (1.2)}\  Assume that  $\delta(f)= 4$ with $d_H(v)= 4$ and
$d_H(u)\geq d_H(w)$. If $n_{7^-}(v)= 1$, then $d_H(u)\geq 10$ by
Claim \ref{4-vertex}, $d_H(w)\geq 6$,  or $4\leq d_H(w)\leq 5$ and
$n_4(w)\geq 1$. Then, $u$ sends $\frac{7}{5}$ to $f$ by (R4) and $w$
sends at least $\frac{4}{5}$ to $f$ by (R1), (R2), (R3) or (R4).
Hence, $w'(f) \geq -3 + 2\times \frac{4}{5}+ \frac{7}{5}= 0$. If
$n_{7^-}(v)= 0$, then $d_H(u), d_H(w)\geq 8$. Then, each of $u, w$
sends $\frac{5}{4}$ to $f$ by (R3) and (R4) and $v$ sends
$\frac{1}{2}$ to $f$ by (R1). Hence, $w'(f) \geq -3 + 2\times
\frac{5}{4}+ \frac{1}{2}= 0$.

{\bf (1.3)}\  Assume that  $\delta(f)= 5$ with $d_H(v)= 5$ and
$d_H(u)\geq d_H(w)$. We further have three possibilities:

$\bullet$ \  $n_5(f)= 3$. Since $G$ contains none of ($A1$),
($A3.1$)
 and ($A4.2$), $n_4(u)= n_4(v)= n_4(w)= 0$ and each vertex in $b(f)$ sends
$\frac{7}{5}$ to $f$ by (R2). Thus, $w'(f)\geq -3 + 3\times
\frac{7}{5}= \frac{6}{5}$.

$\bullet$ \  $n_5(f)= 2$ and $d_H(w)= 5$.

 If $n_4(v)= 1$ or $n_4(w)=
1$, then $d(u)\geq 10$ since $G$ contains no ($A4.2$) and
$d_H(u)\geq 10$ by Claim \ref{9-vertex}. Hence, each of $v$ and $w$
sends $\frac{4}{5}$ or $\frac{11}{12}$ to $f$ by (R2), $u$ sends
$\frac{7}{5}$ to $f$ by (R4) and therefore  $w'(f) \geq -3 + 2\times
\frac{4}{5}+ \frac{7}{5}= 0$. Otherwise, $n_4(v)= n_4(w)= 0$. By
(R2), (R3) or (R4), if $d_H(u)= 6$, then $w'(f) \geq -3 + (2\times
\frac{6}{5}+ 1)= \frac{2}{5}$;  if $d_H(u)= 7$, then $w'(f) \geq -3
+ (2\times \frac{13}{14}+ \frac{8}{7})= 0$; if $d_H(u)= 8$, then
$w'(f) \geq -3 + (2\times 1+ \frac{5}{4})= \frac{1}{4}$;  if
$d_H(u)\geq 9$, then  $w'(f) \geq -3 + (2\times \frac{11}{12}+
\frac{4}{3})=\frac{1}{6}$.

$\bullet$ \    $n_5(f)= 1$ and $d_H(u)\geq d_H(w)\geq 6$.

If $n_4(v)= 1$, then $v$ sends $\frac{4}{5}$ to $f$ by (R2). If
$d_H(w)= 6$, then $d(u)\geq 9$ and $d_H(u)\geq 9$ since $G$ contains
no ($A4.2$). Thus, $w$ sends $1$ to $f$ and $u$ sends $\frac{4}{3}$
to $f$ by (R3) or $u$ sends $\frac{5}{4}$ to $f$ by (R4). Hence,
$w'(f) \geq -3 + \min\{\frac{4}{5} + 1 + \frac{5}{4}, \frac{4}{5} +
1 + \frac{4}{3}\}= \frac{1}{20}$. Otherwise, $ d_H(u)\geq d_H(w)\geq
7$ and each of $w$ and $u$ sends at least $\frac{11}{10}$ by (R3) or
(R4). Hence, $w'(f) \geq -3 + \frac{4}{5} + 2\times \frac{11}{10}=
0$. Now assume that $n_4(v)= 0$. By (R2), (R3) or (R4), if $d_H(u)=
d_H(w)= 6$, then $w'(f) \geq -3 + (1 + 2\times 1)= 0$; if $d_H(u)=
7$ and $d_H(w)= 6$, then $w'(f) \geq -3 + (\frac{6}{7} + 1 +
\frac{8}{7})= 0$; if $d_H(u)= 8$ and $d_H(w)= 6$, then $w'(f) \geq
-3 + (\frac{3}{4} + 1 + \frac{5}{4})= 0$; if $d_H(u)\geq 9$ and
$d_H(w)= 6$, then $w'(f) \geq -3 + (\frac{2}{3} + 1 + \frac{4}{3})=
0$; if $d_H(u)= d_H(w)= 7$, then $w'(f) \geq -3 + (\frac{5}{7} +
2\times \frac{8}{7})= 0$; if $d_H(u)\geq 8$ and $d_H(w)= 7$, then
$w'(f) \geq -3 + (\frac{9}{14} + \frac{8}{7} + \frac{5}{4})=
\frac{1}{28}$; if $d_H(u)\geq 8$ and $d_H(w)= 8$, then $w'(f) \geq
-3 + (\frac{1}{2} + 2\times \frac{5}{4})= 0$; if $d_H(u)\geq
d_H(w)\geq 9$, then $w'(f) \geq -3 + (\frac{1}{3} +
2\times\frac{4}{3})= 0$.

{\bf (1.4)}\ Assume that $\delta(f)\geq 6$. Then  each vertex in
$b(f)$ sends at least $1$ to $f$ by (R3) or (R4). Hence, $w'(f) \geq
-3 + 3\times 1= 0$.

\medskip

\noindent{\bf Case 2.}\  $d(f) = 4$.

Let  $f = [uvst]$. Then $w(f) = d(f)- 6 = -2$. If $n_{6^+}(f)\geq
2$, then by (R3) or (R4), each $6^+$-vertex in $b(f)$ gives at least
$1$ to $f$. Hence, $w'(f) \geq -2 + 2\times 1 = 0$. Now we assume
that $n_{6^+}(f)\leq 1$. Note that if $\delta(f) = 3$, then
$n_{9^+}(f) \geq 2$ and if there is a $4$-vertex $v$ in $b(f)$ with
$n_{7^-}(v)= 0$, then $n_{8^+}(f)\geq 2$. Thus it suffice to suppose
that $\delta(f)\geq 5$,  or $\delta(f)= 4$  and each $4$-vertex $v$
in $b(f)$ having $n_{7^-}(v)= 1$. Assume that $\delta(f)= 4$ with
$d_H(v)= 4$ and $d_H(u)\leq 5\leq 10\leq d_H(s)$. It is easy to see
that $n_4(u)\geq 1$. By (R1), (R2) and (R4), $v$ sends $\frac{4}{5}$
to $f$, $u$ sends $\frac{4}{5}$ to $f$, and $s$ sends $1$ to $f$.
Hence, $w'(f) \geq -2 + 1 + 2\times \frac{4}{5}= \frac{3}{5}$.
Assume that $\delta(f)= 5$ with $d(u)= d(v)= d(s)= 5$ and
$d_H(t)\geq 5$. By (R2), (R3) or (R4), each of $u, v, s, t$ gives at
least $\frac{1}{2}$ to $f$ and $w'(f) \geq -2 + 4\times \frac{1}{2}=
0$.
\medskip

\noindent{\bf Case 3.}\  $d(f) = 5$.

We see that $w(f) = d(f)- 6 = -1$. If $n_{6^+}(f)\geq 1$, then by
(R3) or (R4), each $6^+$-vertex in  $b(f)$ gives at least $1$ to
$f$. Hence, $w'(f) = -1 + 1\times 1 = 0$. So assume that
$n_{6^+}(f)= 0$. This implies that $\delta(f)\geq 5$ and $n_{5}(f)=
5$. Then by (R2), each vertices in  $b(f)$ gives at least
$\frac{1}{2}$ to $f$ and  $w'(f) \geq -1 + 5\times \frac{1}{2}=
\frac{3}{2}$.
\medskip

\noindent{\bf Case 4.}\  $d(f) \ge 6$.

It is obvious that $w(f) = d(f)- 6 \geq 0$.

Let $v\in V(H)$. By Claim \ref{minimum}, we may assume that
$d_H(v)\geq 3$.

If $d_H(v) = 3$, then $w(v) = 2d_H(v)- 6 = 0$.

Assume that $d_H(v) = 4$, then $w(v) = 2d_H(v)- 6 = 2$. If
$n_{7^-}(v)= 1$, then $w'(v)\geq 2- 2\times\frac{4}{5}- 2\times
\frac{1}{5} = 0$ by (R1). Otherwise, $n_{7^-}(v)= 0$ and $w'(v)\geq
2- 4\times\frac{1}{2}= 0$.

Assume that  $d_H(v) = 5$, then $w(v) = 2d_H(v)- 6 = 4$. If
$n_{4}(v)\geq 1$, then $w'(v)\geq 4- 5\times\frac{4}{5} = 0$ by
(R2). Otherwise, $n_{4}(v)= 0$. Let $w_1,w_2,\ldots,w_5$ be the
neighbors of $v$ in a cyclic order.  We need to consider three
subcases:

\medskip

\noindent  {\bf Case 1.}\  $n_{5}(v)\geq 1$ with $d_H(w_1)= 5$.

 Since $G$ contains no ($A_{4.2}$), $n_{6}(v)\leq 2$. This leads to
 the following three possibilities:

{\bf (1.1)}\  If $n_{5}(v)= 2$, then $n_{9^+}(v) = 3$ and
$n'_{9^+}(v) = 3$ as  $G$ contains no ($A_{4.2}$). By symmetry, we
may assume that $d_H(w_2)= 5$ or $d_H(w_3)= 5$. Thus, by (R2),
$w'(f) \geq 4 -\max\{(\frac{7}{5} + 2\times \frac{11}{12} + 2\times
\frac{1}{3}), (4\times \frac{11}{12} + \frac{1}{3})\} = 4-
\max\{3\frac{9}{10}, 4\}= 0$.

{\bf (1.2)}\  If $n_{5}(v) = n_{6}(v)= 1$, then $n_{8^+}(v) = 3$ and
$n'_{8^+}(v) = 3$ as  $G$ contains no ($A_{4.2}$). By symmetry, we
may assume that $d_H(w_2)= 6$ or $d_H(w_3)= 6$. Thus, by (R2),
$w'(f) \geq 4 -\max\{(\frac{6}{5} + 1 + \frac{3}{4} + 2\times
\frac{1}{2}), (2\times 1 + 2\times \frac{3}{4} + \frac{1}{2})\}= 4-
\max\{ 3\frac{19}{20}, 4\}= 0$.

{\bf (1.3)}\  Assume that $n_{5}(v) = 1$ and  $n_{7^+}(v)= 4$. If
$d_H(w_2) =
 d_H(w_5)= 7$, then by (R2), $w'(f)\geq 4 - (2\times \frac{13}{14} + 3\times \frac{5}{7}) =
 0$. If one of $w_2$ and $w_5$ is of degree $7$ in $H$, then by (R2), $w'(f)\geq 4 - (\frac{13}{14} + 1 + 2\times \frac{5}{7} + \frac{9}{14}) =
 0$. Otherwise, $\min\{d_H(w_2), d_H(w_5)\}\geq 8$, then by (R2), $w'(f)\geq 4 - (2\times 1 + 2\times \frac{9}{14} + \frac{5}{7}) = 0$.

\medskip

 \noindent {\bf Case 2.}\  $n_{5}(v)= 0$ and $n_{6}(v)\geq 1$ with $d_H(w_1)= 6$.

 Since $G$ contains no ($A_{4.2}$), $n_{6}(v)\leq 2$. If $n_{6}(v)= 1$ and $n_{7^+}(v)= 4$, then by (R2),
$w'(f)\geq 4 - (2\times \frac{6}{7} + 3\times \frac{5}{7}) =
\frac{1}{7}$. Otherwise, $n_{6}(v)= 2$ and $n_{7^+}(v)=3$.  Since
$G$ contains no ($A_{4.2}$), $n_{8^+}(v)\geq 1$, or $n_{7}(v)= 3$
with $uv_5\not\in E(H)$ and $d_H(w_2)= 6$ or $d_H(w_3)= 6$. If
$n_{7}(v)= 3$, then by (R2), $v$ sends at most $\frac{1}{2}$ to the
face  whose boundary contains $w_5, v, w_1$. It follows that
 $w'(f)\geq 4 -(\frac{1}{2}
+ 1 + \frac{6}{7} + 2\times\frac{5}{7}) = \frac{3}{14}$ if
$d_H(w_2)= 6$,   or $w'(f)\geq 4 -(\frac{1}{2} + 3\times \frac{6}{7}
+ \frac{5}{7}) = \frac{3}{14}$ if $d_H(w_3)= 6$. Otherwise,
$n_{8^+}(v)\geq 1$ and we need to consider the following subcases by
symmetry.

$\bullet$ \  If $d_H(w_2)= 6$ and $d_H(w_3)\geq 8$, then by (R2),
$w'(f)\geq 4 - (1 + \frac{6}{7} + \frac{5}{7} + \frac{3}{4} +
\frac{9}{14}) = \frac{1}{28}$.

$\bullet$ \   If $d_H(w_2)= 6$ and $d_H(w_4)\geq 8$, then by (R2),
$w'(f)\geq 4 - (1 + 2\times \frac{6}{7} + 2\times \frac{9}{14}) =
0$.

$\bullet$ \  If $d_H(w_3)= 6$ and $d_H(w_2)\geq 8$, then by (R2),
$w'(f)\geq 4 - (2\times \frac{6}{7} +2\times \frac{3}{4} +
\frac{5}{7}) = \frac{1}{14}$.

$\bullet$ \  If $d_H(w_3)= 6$ and $d_H(w_4)\geq 8$, then by (R2),
$w'(f)\geq 4 - (3\times \frac{6}{7} + \frac{3}{4} + \frac{9}{14}) =
\frac{1}{28}$.

\medskip

\noindent {\bf Case 3.}\  $n_{5}(v)= n_{6}(v)= 0$ and $n_{7^+}(v)=
5$. Then, by (R2), $w'(f)\geq 4 - 5\times \frac{5}{7}= \frac{3}{7}$.

If $6\leq d_H(v)\leq 9$, then by (R3), $w'(f)= 2d_H(v)- 6-
d_H(v)\times \frac{2d_H(v)- 6}{d_H(v)}= 0$.

If $d_H(v)\geq 12$, then by (R4), $w'(f)\geq 2d_H(v)- 6-
d_H(v)\times \frac{3}{2}= \frac{1}{2}d_H(v)- 6\geq \frac{1}{2}\times
12- 6= 0$.

Assume that  $d_H(v)= 10$, then $w(f)= 2d_H(v)- 6= 14$. Let $t$ the
number of faces $f$ incident to $v$ with  $d(f)\geq 4$, or $d(f)= 3$
and $\delta(f)\geq 6$. If $t\geq 2$, then by (R4), $w'(f)\geq 14-
8\times \frac{3}{2}- 2\times 1= 0$. Otherwise, $t\leq 1$. Note that
if $n_{5^-}(v)\leq 4$, then $t\geq 2$. Recall that if $d(v)\neq 10$,
then $n_{5^-}(v)\leq 2\leq 4$ by Claim \ref{10+345}. Hence, it
suffice to assume that $d(v)= 10$ with $v_1, v_2, \ldots, v_{10}$ as
the neighbors of $v$ in a clockwise order, $n_{5^-}(v)\geq 5$ and
$m_{4^+}(v)\leq t\leq 1$. If $v$ is not incident to a bad $3$-face,
then $w'(f)= 14- 10\times \frac{7}{5}= 0$ by (R4). Otherwise,
$n_3(v)\geq 1$ and $v$ is incident to some bad $3$-faces. Since $G$
contains no ($A_{3.3}$), there is no $3$-vertex $w$ such that the
faces $f_1, f_2$ which is incident to $v, w$, are both $3$-faces. By
the previous discussion, we may assume that $d_H(v_1)= 3$,
$v_1v_2\in E(H)$, $f_{10}=[\cdots v_{10}vv_1\cdots]$ is a $4^+$-face
and the other faces incident to $v$ are all of degree $3$. Let
$m'_3(v)$ be the number of bad $3$-faces incident to $v$. Clearly,
$1\leq m'_3(v)\leq 2$,  and $m'_3(v)= 2$ if and only if
$d_H(v_{10})= 3$. If $m'_3(v)= 2$, then $w'(f)\geq 14-
(2\times\frac{3}{2} + 1 + 7\times \frac{7}{5})=\frac{1}{5}$ by (R4).
Otherwise, $m'_3(v)= 1$ and hence $w'(f)\geq 14- (1\times\frac{3}{2}
+ 1 + 8\times \frac{7}{5})=\frac{3}{10}$.

Finally assume that $d_H(v)= 11$, then $w(f)= 2d_H(v)- 6= 16$. If
$n_{3}(v)\leq 3$, then the number of bad $3$-faces incident to $v$
is at most $6$. Hence, by (R4), $w'(f)\geq 16- 6\times \frac{3}{2}-
5\times \frac{7}{5}= 0$. Note that if $d(v)\neq 11$, then
$n_{3}(v)\leq n_{5^-}(v)\leq 3$ by Claim \ref{10+345}. Now we assume
that $d(v)= 11$ with $w_1, w_2, \ldots, w_{11}$ as the neighbors of
$v$ in a clockwise order, $n_{3}(v)\geq 4$ and $v$ is incident to at
least $7$ bad $3$-faces. If $v$ is incident to a $4^+$-face,  or to
a $3$-face $f$ with $\delta(f)\geq 6$, then $w'(f)= 16- 10\times
\frac{3}{2}- 1= 0$. Otherwise, $m_3(v)= 11$ and $\delta(f)\leq 5$
for each $f$ incident to $v$, i.e., for $1\leq i\leq 11$ (where
$w_{12}= w_1$), $\min \{d_H(w_i), d_H(w_{i + 1})\}\leq 5$.

Note that if $\min \{d_H(w_i), d_H(w_{i + 1})\}= 3$, then
$\max\{d_H(w_i), d_H(w_{i + 1})\}\geq 10$ by Claim \ref{3-face} for
each $1\leq i\leq 11$ (where $w_{12}= w_1$). Thus, $n_3(v)\leq 5$.
If $n_3(v)= 5$, then we may assume that $d_H(w_i)= 3$ for $i= 3, 5,
7, 9, 11$ and $d_H(w_j)\geq 10$ for $j= 1, 2, 4, 6, 8, 10$.
Obviously, $f_1=[w_1vw_2]$ is of degree $3$ with $\delta(f)\geq 10$,
a contradiction. Otherwise, $n_3(v)= 4$. Further, for $1\leq i \leq
11$ (where $w_{12}= w_1$, $w_{13}= w_2$, $w_{14}= w_3$), if
$d_H(w_i)= d_H(w_{i + 3}) = 3$, then $d_H(w_{i + 1})\geq 10,
d_H(w_{i + 2})\geq 10$ by Claim \ref{3-face} and $f_{i + 1}= [w_{i +
1}vw_{i + 2}]$ is of degree $3$ with $\delta(f_{i + 1})\geq 10$, a
contradiction. Hence, it suffice to assume that $d_H(w_i)= 3$ for
$i= 4, 6, 8, 10$; $d_H(w_j)\geq 10$ for  $j= 3, 5, 7, 9, 11$ and $4
\leq d_H(w_1), d_H(w_2)\leq 5$.

Since $G$ contains no ($A_{4.1}$), there is no $4$-vertex adjacent
to a $4$-vertex and $11$-vertex at the same time and
$\max\{d_H(w_1), d_H(w_2)\}= 5$. Without loss of generality, assume
that $d_H(w_1)= 5$. If $d_H(w_2)= 4$, then $n_4(w_1)\geq 1$ and
$n_{7^-}(w_2)= 1$. By (R4), $w'(f)= 16 - 8\times \frac{3}{2}-
2\times \frac{7}{5}- \frac{11}{10}= \frac{1}{10}$. Otherwise,
$d_H(w_2)= 5$. If $n_4(w_1)\geq 1$ and $n_4(w_2)\geq 1$, then By
(R4), $w'(f)= 16 - 8\times \frac{3}{2}- 2\times \frac{11}{10}-
\frac{7}{5}= \frac{2}{5}$. If $n_4(w_1)= n_4(w_2)= 0$, then By (R4),
$w'(f)= 16 - 8\times \frac{3}{2}- 3\times \frac{4}{3}= 0$.
Otherwise, by (R4), $w'(f)= 16 - 8\times \frac{3}{2}- \frac{7}{5}-
\frac{4}{3}- \frac{11}{10}= \frac{1}{6}$. \qed

%%%%%%%%%%%%%%%%%%%%%%%%%%%%%%%%%%%%%%%%%%%%%%%%%%%%%%%%%%%%%%
\section{Acyclic chromatic indices}
%%%%%%%%%%%%%%%%%%%%%%%%%%%%%%%%%%%%%%%%%%%%%%%%%%%%%%%%%%%%%%

In this section, we discuss the acyclic chromatic indices of planar
graphs. The proof of the main result requires the following  lemmas.

\begin{lemma}\label{MLSC09} {\rm (\cite{MLSC09})}
If $G$ is a graph with $\Delta\leq4$, then $a'(G)\le 7$.

\end{lemma}

Assume that $c$ is a partial acyclic edge $k$-coloring of a graph
$G$ using the color set $C = \{1, 2, \ldots, k\}$. For a vertex
$v\in V(G)$, we use $C(v)$ to denote the set of colors assigned to
edges incident to $v$ under $c$. If the edges of a cycle are
alternatively colored with colors $i$ and $j$, then we call such
cycle an {\em $(i, j)$-cycle}. If the edges of a path $ux\ldots v$
are alternatively colored with colors $i$ and $j$, then we call such
path an {\em $(i, j)_{(u, v)}$-path}.

\begin{lemma}\label{exchange}

If $P= uv_1v_2\ldots v_kv_{k + 1}$ is an $(i, j)_{(u, v_{k +
1})}$-path in $G$ with $c(uv_1)= i$, $j\not\in C(u)$ and $w\not\in
V(P)$, then there is no $(i, j)_{(u, w)}$-path in $G$.

\end{lemma}

\proof\ Otherwise, assume that there is an $(i, j)_{(u, w)}$-path
$Q=uw_1w_2\ldots w_m$ in $G$, where $w=w_m$ and $m\ge 1$. Note that
some $w_i$ may be identical to some $v_j$, $1\le i\le k+1$ and $1\le
j\le m-1$. Since $j\notin C(u)$ and $w\not\in V(P)$, it is easy to
see that there exist  a vertex $v_i$ and  some vertex $w_i$ such
that $v_iv_{i+1}$ and $v_iw_i$ have same color $i$ or $j$, which
contradicts the fact that the edge coloring considered is proper.
  \qed

\begin{theorem}\label{thplanar}
If $G$ is a planar graph, then $a'(G)\leq\Delta + 7$.
\end{theorem}
\proof\ The proof is proceeded by induction on the edge number
$|E(G)|$. If $|E(G)|\leq \Delta + 7$, then $G$ is obviously
acyclically edge $(\Delta + 7)$-colorable. Assume that $G$ is a
planar graph   with $|E(G)|\geq \Delta + 8$. If $\Delta\leq 4$, then
$a'(G)\leq\Delta + 3$ by Lemma  \ref{MLSC09}. Thus, we may assume
that $\Delta\geq 5$ and $G$ is $2$-connected. By Lemma \ref{planar},
$G$ contains one of the configurations ($A_{1}$)-($A_{4}$). In the
following, we deal with each of the configurations
($A_{1}$)-($A_{4}$).

\medskip
\noindent{\bf (A$_{1}$)}\ There is a path $uvw$ such that $d(v)= 2$
and $d(u)\le 9$.

Let $v,u_1,u_2,\ldots,u_{d(u)-1}$ be the neighbors of $u$ and let $H
= G - uv$. Then $H$ is a planar graph and $\Delta(H) \ge
\Delta-1\geq 4$. By the induction assumption or Lemma  \ref{MLSC09},
$H$ has an acyclic edge $(\Delta + 7)$-coloring $c$ using the color
set $C=\{1, 2, \ldots, \Delta + 7\}$.

If $c(vw)\not\in C(u)$, then we color $uv$ with a color in
$C\backslash (C(u)\cup C(v))$. Otherwise, $c(vw)\in C(u)$ and assume
that $c(vw)= c(uu_1)= 1$, $c(uu_i)= i$, $i= 2, 3, \ldots, d(u)- 1$.
If $H$ contains no $(1, a)_{(w, u_{1})}$-path for some color $a\in
C\backslash (C(u)\cup C(v))$, we color $uv$ with $a$. Otherwise,
assume that $H$ contains a $(1, i)_{(w, u_{1})}$-path for any $i\in
C\backslash (C(u)\cup C(v))$. This implies that $d(u) = 9$ and
$C(w)= C(u_{1})= (C\backslash C(v))\cup \{1\}= C\backslash\{2, 3, 4,
5, 6, 7, 8\}$. First, we recolor $vw$ with $2$. Similarly to the
previous argument, we may assume that $H$ contains a $(2,
i)_{(u_{2}, w)}$-path for any $i\in C\backslash (C(u)\cup C(v))$ and
$C(u_{2})= C\backslash\{1, 3, 4, 5, 6, 7, 8\}$. So we switch the
colors of $uu_{1}$ and $uu_{2}$, then color $uv$ with $9$. Clearly,
no bichromatic cycles are produced in $G$ under the constructed
coloring by the previous assumption.

\medskip

\noindent {\bf (A$_{2}$)}\ There is a vertex $u$ with $n_{2}(u)\geq
1$ and $n_{8^{-}}(u)\geq d(u)-8$. Suppose that $u_{1}, u_{2},
\ldots, $ $u_{d(u)-1}$, $v$ are the neighbors of $u$ such that
$d(u_{1})\geq d(u_{2})\geq \cdots \geq d(u_{d(u)-1})\geq d(v)= 2$.
Let $w$ be the neighbor of $v$ different from $u$. For $1\le i\le
d(u)-1$, if $d(u_i)=2$, then we use $x_i$ to denote the neighbor of
$u_i$ different from $u$. Then at least one of the configurations
($A_{2.1}$)-($A_{2.2}$) holds.
\medskip

Without loss of generality, assume that $c(uu_i)= i$ for $1\leq
i\leq d(u)-1$ and let $S_{i^-}= \{c(uu_j)| d(u_j)\leq i\}$. Since
$d(u)\leq\Delta$, $|C\backslash C(u)|\geq 8$. If $c(vw)\not\in
C(u)$, then we color $uv$ with a color in $C\backslash (C(u)\cup
C(v))$. Otherwise, assume that $c(vw)\in C(u)$ and $H$ contains a
$(c(vw), i)_{(u, w)}$-path for any $i\in C\backslash C(u)$,
$(C\backslash C(u))\subseteq C(w)$. If $c(vw)= c(uu_i)\in S_{8^-}$,
we color $uv$ with a color in $(C\backslash C(u))\backslash C(u_i)$.
Otherwise, $c(vw)\in C(u)\backslash S_{8^-}$. Further, we may assume
that $S_{8^-}\subseteq C(w)$. Otherwise, we recolor $vw$ with a
color in $(S_{8^-})\backslash C(w)$ and reduce the proof to the
 case in which $c(vw)\in S_{8^-}$. So we assume that:

\medskip

\noindent($*_{2.1}$)\ $c(vw)\in C(u)\backslash S_{8^-}$ and
$S_{8^-}\cup (C\backslash C(u))\subseteq C(w)\backslash \{c(vw)\}$.

\medskip
 \noindent{\bf (A$_{2.1}$)}\ $n_{8^{-}}(u)\geq d(u) - 7$.
\medskip

By ($*_{2.1}$), $(C\backslash \{c(uu_1), c(uu_2), \cdots,
c(uu_7)\})\subseteq C(w)\backslash \{c(vw)\}$ and $d(w)=
|C(w)|\geq\Delta + 1$.

\medskip
\noindent{\bf (A$_{2.2}$)}\ $n_{8^{-}}(u) = d(u) - 8$, and
$n_{2}(u)\geq d(u)-9$.

\medskip

Without loss of generality, we assume that $d(u_9)= 8$.  By
($*_{2.1}$), we assume that $d(u_i)\geq 9$, $i= 1, 2, \cdots, 8$,
$c(vw) = 1$ and $C(w) = C\backslash\{2, 3, \ldots, 8\}$. Further, we
assume that $H$ contains a $(1, i)_{(u_1, w)}$-path for any $ i\in
C\backslash C(u)$. If $H$ contains no $(i, j)_{(u_i, w)}$-path for
some $j\in C\backslash C(u)$ and some $i\in \{2, 3, \cdots, 8\}$, we
recolor $vw$ with $i$ and color $uv$ with $j$. Otherwise, we may
assume that:

\medskip

\noindent($*_{2.2}$)\ For any $i\in \{1, 2, \cdots, 8\}$ and $ j\in
C\backslash C(u)$, $H$ contains an $(i, j)_{(u_i, w)}$-path and
$(C\backslash C(u))\subseteq C(u_i)$.
\medskip

Assume that $H$ contains no $(i, 9)_{(u_i, w)}$-path for some color
$i\in \{1, 2, \ldots, 8\}$. Let $C(u_9)= \{9, a_1, a_2, \ldots,
a_7\}$. If $\{a_1, a_2, \ldots, a_7\}\cap S_2= \emptyset$, then we
recolor $uu_9$ with a color in $C\backslash(C(u)\cup C(u_9))$.
Otherwise, assume that $\{a_1, a_2, \ldots, a_7\}\cap S_2 \neq
\emptyset$. Let $I= \{i | i \in \{10, \ldots, d(u)- 1\}$ such that
$c(uu_i)\in \{a_1, a_2, \ldots, a_7\} \}$. Then we recolor $uu_9$
with a color in $C\backslash(C(u)\cup C(u_9)\cup \{C(u_i)| i\in
I\})$,  $vw$ with $i$ if $i\neq 1$ and color $uv$ with $9$. So we
may assume that:

\medskip

\noindent($*_{2.3}$)\ $H$ contains a $(i, 9)_{(u_i, w)}$-path for
any $ i\in \{1, 2, \ldots, 8\}$.

\medskip

Assume that $H$ contains no $(i, j)_{(u_i, w)}$-path for some $i\in
\{1, 2, \ldots, 8\}$ and some $j= c(uu_j)\in S_2$. Let $C(u_j)= \{j,
a\}$. If $a\not\in S_2\cup \{9\}$, we recolor $uu_j$ with a color in
$C\backslash(C(u)\cup C(u_j))$; If $a=\{c(uu_k)\}\in S_2\cup \{9\}$,
we recolor $uu_j$ with a color in $C\backslash(C(u)\cup C(u_j)\cup
C(u_k))$. Recolor $vw$ with $i$ if $i\neq 1$ and color $uv$ with
$j$. So we may assume that:

\medskip

\noindent($*_{2.4}$)\ $H$ contains an $(i, j)_{(u_i, w)}$-path for
any $i\in \{1, 2, \ldots, 8\}$ and any $ j\in S_2$.

\medskip

By ($*_{2.2}$), ($*_{2.3}$), and ($*_{2.4}$), we may assume that:

\medskip

\noindent($*_{2.5}$)\ For any $j\in C\backslash\{1, 2, \ldots, 8\}$,
$H$ contains an $(i, j)_{(u_i, w)}$-path, and $C(u_i)\backslash$
$\{i\} = C(w)$ $\backslash \{c(vw)\}$ $= C\backslash$ $\{1, 2,
\ldots, 8\}$, $i = 1, 2, \ldots, 8$.

\medskip

\begin{claim}\label{abc}
There exist $ i_0, j_0 \in \{1, 2, \ldots, 8\}$ such that $H$
contains neither an $(i_0, 9)_{(u_{j_0}, u_9)}$-path nor a $(j_0,
9)_{(u_{i_0}, u_9)}$-path.
\end{claim}

\proof Since $d(u_9) = 8$, $\{1, 2, \ldots, 8\}\backslash
C(u_9)\neq\emptyset$. Assume that $1\not\in C(u_9)$. If $H$ contains
no $(2, 9)_{(u_1, u_9)}$-path, $\{i_0, j_0\}= \{1, 2\}$. If $H$
contains no $(3, 9)_{(u_1, u_9)}$-path, $\{i_0, j_0\}= \{1, 3\}$.
Otherwise, $H$ contains a $(2, 9)_{(u_1, u_9)}$-path and a $(3,
9)_{(u_1, u_9)}$-path. It follows that $H$ contains neither a $(2,
9)_{(u_3, u_9)}$-path nor a $(3, 9)_{(u_2, u_9)}$-path, thus $\{i_0,
j_0\}= \{2, 3\}$. \qed

By symmetry, assume that $\{i_0, j_0\} = \{1, 2\}$. We switch the
colors of $uu_{1}$ and $uu_{2}$. If no bichromatic cycles are
produced, we color $uv$ with some color in $C\backslash C(u)$ and we
are done. Otherwise, by the previous assumption, it is easy to
conclude that the bichromatic cycle in $G$ must be a $(1, i)_{(u_2,
u_i)}$-cycle or a $(2, j)_{(u_1, u_j)}$-cycle for some $i, j\in
S_2$. Let

\begin{center}
$T_1= \{u_i |$ $G$ contains a $(1, i)_{(u_2, u_i)}$-cycle$\}$,
\end{center}

\begin{center}
$T_2= \{u_i |$ $G$ contains a $(2, i)_{(u_1, u_i)}$-cycle$\}$.
\end{center}

If $|T_1| \geq 2k$, $k\geq 1$, assuming that $u_{i_1}, u_{i_2}\in
T_1$, then we switch the colors of $uu_{i_1}$ and $uu_{i_2}$, so
that $G$ contains neither $(1, i_1)_{(u_2, u_{i_1})}$-cycle nor $(1,
i_2)_{(u_2, u_{i_2})}$-cycle, and no other new bichromatic cycles
are produced. Since $\min\{|T_1|, |T_2|\}\geq 1$, we assume that
$|T_1|= 1$, $|T_2|\leq 1$. Without loss of generality, assume that
$u_{10}\in T_1$, $c(u_{10}x_{10}) = 1$, and $10\in C(x_{10})$.

Since $d(x_{10})\leq \Delta$, $(C\backslash\{2\})\backslash
C(x_{10})\neq \emptyset$, assume that $p\in
(C\backslash\{2\})\backslash C(x_{10})$ and recolor $u_{10}x_{10}$
with $p$. If $p\in (C\backslash C(u))\cup \{3, 4, \ldots, 8\}$, then
let $S'_2 = S_2$ and $C'(u)= C(u)$. If $p= c(uu_i)\in S_2\cup
\{9\}$, recolor $uu_{10}$ with a color $j_1\in (C\backslash
C(u))\backslash C(u_i)$ and let $S'_2= (S_2\backslash \{10\})\cup
\{j_1\}$, $C'(u)=(C(u)\backslash\{10\})\cup \{j_1\}$. If $|T_2| =
0$, then we color $uv$ with $10$.  Otherwise, $|T_2| =1$, $u_{11}\in
T_2$ and $c(u_{11}x_{11}) = 11$. Since $d(x_{11})\leq \Delta$,
$(C\backslash\{1\})\backslash C(x_{11})\neq \emptyset$, assume that
$q\in (C\backslash\{1\})\backslash C(x_5)$ and recolor
$u_{11}x_{11}$ with $q$. If $q\in (C\backslash C'(u))\cup \{3, 4,
\ldots, 8\}$, then let $C''(u) = C'(u)$. If $q= c(uu_i)\in S'_2\cup
\{9\}$, we recolor $uu_{11}$ with a color $j_2\in (C\backslash
C'(u))\backslash C(u_i)$ and $C''(u) = (C'(u)\backslash\{11\})\cup
\{j_2\}$. Then, we color $uv$ with a color in $C\backslash C''(u)$.

\medskip

\noindent{\bf (A$_{3}$)}\ There is a $3$-vertex $u$ adjacent to a
vertex $v$ such that one of ($A_{3.1}$)-($A_{3.3}$) holds:

Let $v_1, v_2, \ldots, v_{d(v)}$ be the neighbors of $v$, where
$v_{8}= u_2$ in ($A_{3.2}$), ($v_{8}= u_1$, $v_{9}= u_2$ in
($A_{3.3}$)) and $v_{d(v)}= u$. Let $H = G - uv$. Then $H$ is a
planar graph and $\Delta(H) \ge \Delta-1\geq 4$. By the induction
assumption or Lemma \ref{MLSC09}, $H$ has an acyclic edge $(\Delta +
7)$-coloring $c$ using the color set $C=\{1, 2, \ldots, \Delta +
7\}$. Note that $|C\backslash C(v)|\geq 8$. If $C(u)\cap C(v)=
\emptyset$, then color $uv$ with a color in $C\backslash (C(u)\cup
C(v))$. Otherwise, $C(u)\cap C(v)\neq \emptyset$, i.e., $1\leq|
C(u)\cap C(v)|\leq 2$ and $c(uu_i)= i$, $i= 1, 2$.

If, for any $i\in C(u)\cap C(v)$, $H$ contains no $(i, j)_{(u,
v)}$-path for some $j\in C\backslash (C(u)\cup C(v))$, then we can
color $uv$ with $j$. Otherwise, we assume that:

\medskip

\noindent($*_{3.1}$)\ For any $j\in C\backslash (C(u)\cup C(v))$,
$H$ contains an $(i, j)_{(u, v)}$-path for some   $i\in C(u)\cap
C(v)$.

Next we need to consider the following two subcases.

\medskip
\noindent {\bf (I)}\ $|C(u)\cap C(v)|=1$.
\medskip

\noindent{\bf (A$_{3.1}$)} \ $d(v)\leq 8$.
\medskip

Assume that $C(v) = \{1, 2, \ldots, d(v)\}\backslash \{2\}$, and
$c(uu_1) = c(vv_1) = 1$. By ($*_{3.1}$), $d(v)= 8$, $C(u_1) = C(v_1)
= C\backslash \{2, 3, \ldots, 8\}$. We recolor $uu_1$ with $3$,
$vv_1$ with $2$, and color $uv$ with $1$. Clearly, no bichromatic
cycles are produced in $G$ under the constructed coloring.

\medskip

\noindent{\bf (A$_{3.2}$)} \ $d(v)= 9$ and $uu_2, vu_2 \in E(G)$.
\medskip

{\bf (3.2.1)}\ $C(v) = \{1, 2, \ldots, d(v)\}\backslash \{2\}$.

Assume that $c(vu_2)= 1$. By ($*_{3.1}$), $d(u_2)= \Delta$,
$C(v)\cup C(u_2)= C$, $C(v)\cap C(u_2)= \emptyset$ and
$C(u_2)\backslash \{2\}\subseteq C(u_1)$. Since $d(u_1)\leq\Delta$,
there exists $a_1\in C(v)\backslash C(u_1)$. We recolor $uu_2$ with
$a_1$ and color $uv$ with $2$.

Assume that $c(vv_i)= i$, $i\in \{1, 2, \ldots, 7\}\backslash
\{2\}$, and $c(vu_2)= 8$, $c(vv_2)= 9$. By ($*_{3.1}$), $C\backslash
(2, 3, \ldots, 9)\subseteq C(u_1)$. Since $d(u_1)\leq \Delta$,
$|\{2, 3, \ldots, 9\}\cap C(u_1)|\leq 1$. If $8\not\in C(u_1)$, we
recolor $uu_1$ with $8$ and no bichromatic cycle is produced since
$2\not\in C(v)$, then the proof is similar to the case $c(vu_2)=
c(uu_1)= 1$. Otherwise, $C(u_1)= (C\backslash (C(u)\cup C(v)))\cup
\{1, 8\}$. If there exists $a_2\in (C\backslash (C(u)\cup
C(v)))\backslash C(u_2)$, we recolor $uu_2$ with $a_2$ and color
$uv$ with $2$. Otherwise, $C(u_2)= (C\backslash (C(u)\cup C(v)))\cup
\{2, 8\}$. Then, we recolor $uu_2$ with $1$, $uu_1$ with $3$ and
color $uv$ with $2$.

{\bf (3.2.2)}\ $C(v) = \{1, 2, \ldots, d(v)\}\backslash \{1\}$.

Assume that $c(vv_i)= i$, $i\in \{1, 2, \ldots, 7\}\backslash
\{1\}$, and $c(vu_2)= 8$, $c(vv_1)= 9$. By ($*_{3.1}$), $d(u_2)=
\Delta$, $C(v)\cup C(u_2)= C$, $C(v)\cap C(u_2)= \emptyset$. Then,
we recolor $uu_1$ with a color in $C\backslash (C(u_1)\cup \{2,
8\})$ and color $uv$ with $1$.

\medskip

\noindent{\bf ($A_{3.3}$)}\ $d(v)= 10$, $n_{5^-}(v)\geq 5$ and
$uu_1, vu_1, uu_2, vu_2 \in E(G)$.
\medskip

Without loss of generality assume that $C(v) = \{1, 2, \ldots,
d(v)\}\backslash \{2\}$ and $c(vv_i)= i$, $i\in \{3, 4, \ldots,
7\}$, $c(vu_1)= 8$, $c(vv_2)= 10$.

{\bf (3.3.1)}\ $c(vu_2)= 1$ and $c(vv_1)= 9$.

By ($*_{3.1}$), $(C\backslash (C(u)\cup C(v)))\subseteq C(u_1)\cap
C(u_2)$. Since $d(u_1), d(u_2)\leq \Delta$, $|\{3, 4, \ldots,$ $7,
9, 10\}\cap$ $C(u_1)|\leq 1$, and $|\{3, 4, \ldots, 7, 9, 10\}\cap$
$C(u_2)|\leq 1$. Thus, there exists $a_3\in $ $\{3, 4, \ldots, 7, 9,
10\}\backslash$ $(C(u_1)\cup C(u_2))$. Then, we recolor $uu_2$ with
$a_3$ and color $uv$ with $2$.

{\bf (3.3.2)}\ $c(vu_2)= 9$ and $c(vv_1)= 1$.

By ($*_{3.1}$), $(C\backslash (C(u)\cup C(v)))\subseteq C(u_1)\cap
C(v_1)$. If $9\not\in C(u_1)$, we recolor $u_1u$ with $9$ and the
proof is similar to (3.3.1). Otherwise, $9\in C(u_1)$. Since
$d(u_1)\leq \Delta$, $d(u_1)= \Delta$, and $C(u_1)= (C\backslash
(C(u)\cup C(v)))\cup \{1, 8, 9\}$. Obviously, $2\not\in C(u_1)$. If
there exists $a_4\in (C\backslash C(u)\cup C(v))\backslash C(u_2)$,
recolor $uu_2$ with $a_4$ and color $uv$ with $2$. By ($*_{3.1}$)
and Lemma \ref{exchange}, no bichromatic cycles are produced.
Otherwise, $(C\backslash C(u)\cup C(v))\subseteq C(u_2)$ and there
exists $a_5\in \{3, 4, \ldots, 7, 10\}\backslash C(u_2)$. Then, we
recolor $uu_2$ with $a_5$ and color $uv$ with $2$.

\medskip
\noindent {\bf (II)}\ $|C(u)\cap C(v)|=2$.
\medskip

If there exists $a_6\in C\backslash (C(u)\cup C(v)\cup C(u_{i_0}))$
for some $i_0\in \{1, 2\}$, then we recolor $uu_{i_0}$ with $a_6$
and reduce the proof to the Case I. Otherwise, we assume that:

\medskip

\noindent($*_{3.2}$)\ For each $i\in \{1, 2\}$, $C\backslash
(C(u)\cup C(v)\cup C(u_i)) = \emptyset$, i.e., $(C(u)\cup C(v)\cup
C(u_i))= C$.
\medskip

Note that in ($A_{3.1}$) and ($A_{3.2}$), $C\backslash (C(u)\cup
C(v)\cup C(u_2)) \neq \emptyset$. So we assume that $d(v)= 10$,
$n_{5^-}(v)\geq 5$ and $uu_i, vu_i \in E(G)$, $i= 1, 2$. Further,
assume that $C(v)= \{1, 2, \ldots, 9\}$ and $c(vv_i)= i$, $i= \{1,
2, \ldots, 9\}\backslash \{2, 8\}$.

$\bullet$ \ $c(vu_1)= 8$, $c(vv_2)= 2$.

By ($*_{3.2}$), $d(u_i)= \Delta$, $i= 1, 2$ and $C(u_1)=
(C\backslash$ $\{1, 2, \ldots,$ $9\})\cup$ $\{1, 8\}$, $C(u_2)=
(C\backslash \{1, 2, \ldots, 9\})\cup \{2, 9\}$. Since
$n_{5^-}(v)\geq 5$, there exists $j= c(vv_j), k= c(vv_k)$ and
$d(v_j), d(v_k)\leq 5$. Without loss of generality, assume that $j=
1$ and $k= 2$ (Otherwise, we can recolor $uu_1$ with $j$, $uu_2$
with $k$ if $\{j, k\}\neq \{1, 2\}$ and the proof is similar.
Clearly, no bichromatic cycles are produced and $C(u)\cap C(v)= \{j,
k\}$.) By ($*_{3.1}$), for any $ j\in C\backslash (C(u)\cup C(v))$,
$H$ contains an $(i, j)_{(u, v)}$-path for some color $i\in \{1,
2\}$. Hence, $\Delta = 10$, $C(u_1)= \{1, 10, 11, 12, 13\}$ and
$C(u_2)= \{2, 14, 15, 16, 17\}$. Then, switch the colors of $uu_1$
and $uu_2$, and color $uv$ with $10$. By ($*_{3.1}$) and Lemma
\ref{exchange}, no bichromatic cycles are produced.

$\bullet$ \ $c(vu_1)= 2$, $c(vv_2)= 8$.

By ($*_{3.2}$), $d(u_i)= \Delta$, $i= 1, 2$ and $C(u_1)=
(C\backslash$ $\{1, 2, \ldots,$ $9\})\cup \{1, 2\}$, $C(u_2)=
(C\backslash \{1, 2, \ldots, 9\})\cup \{2, 9\}$. We recolor $uu_2$
with $8$ and the proof is similar to the case $c(vu_1)= 8$ and
$c(uu_2)= 2$.

\medskip

\noindent{\bf (A$_{4}$)}\ There is a vertex $v$ adjacent to $u, v_2,
\ldots, v_{d(v)}$ with  $d(u)\leq d(v_2) \leq \ldots \leq
d(v_{d(v)})$ such that one of ($A_{4.1}$)-($A_{4.2}$) holds.
\medskip

Let $u_1, u_2, \ldots, u_{d(u)- 1}\neq v$ be the other neighbors of
$u$. Let $H = G - uv$. Then $H$ is a planar graph and $\Delta(H) \ge
\Delta-1\geq 4$. By the induction assumption or Lemma  \ref{MLSC09},
$H$ has an acyclic edge $(\Delta + 7)$-coloring $c$ using the color
set $C=\{1, 2, \ldots, \Delta + 7\}$ with $c(uu_i)= i$, $i= 1, 2,
\ldots, u_{d(u)- 1}$. Note that $|C\backslash C(u)|\geq 8$,
$|C(v)|\leq d(v)\leq 5$ and $C\backslash (C(u)\cup C(v))\neq
\emptyset$. If $C(u)\cap C(v)= \emptyset$, then we color $uv$ with a
color in $C\backslash (C(u)\cup C(v))$. Otherwise, $C(u)\cap
C(v)\neq \emptyset$.

If, for any $i\in C(u)\cap C(v)$, $H$ contains no $(i, j)_{(u,
v)}$-path for some $j\in C\backslash (C(u)\cup C(v))$, then we can
color $uv$ with $j$. Otherwise, we assume that:

\medskip

\noindent($*_{4.1}$)\ For any $ j\in C\backslash (C(u)\cup C(v))$,
$H$ contains an $(i, j)_{(u, v)}$-path for some   $i\in C(u)\cap
C(v)$.

Next we need to consider the following four subcases.

\medskip
\noindent {\bf Case 1}\ $|C(u)\cap C(v)|=1$ and $1= c(uu_1)\in
C(v)$.
\medskip

By ($*_{4.1}$), we only need to consider two cases:  (i) $d(u)= 7$,
$d(v)= 4$ and $C(v)= \{1, 7, 8\}$ with $c(w_1)= 1$, $w_1, w_2\in
\{v_1, v_2, v_3\}$; (ii) $d(u)= 6$, $d(v)= 5$ and $C(v)= \{1, 6, 7,
8 \}$ with $c(w_1)= 1$, $w_1, w_2, w_3\in \{v_1, v_2, v_3, v_4\}$.
By ($*_{4.1}$), we can assume that $C(u_1) = C(w_1)= (C\backslash
\{1, 2, \ldots, 8\})\cup \{1\}$. Then, we recolor $uu_1$ with $7$,
$vw_1$ with $2$ and color $uv$ with $1$. Since $7, 8\not\in C(w_1)$
and $1\not\in C(u_1)$, no bichromatic cycles are produced.

Otherwise, $|C(u)\cap C(v)|\geq 2$ and $|C\backslash (C(u)\cup
C(v))|= \Delta + 7- |C(u)\cup C(v)|= \Delta + 7 -(|C(u)| + |C(v)| -
|C(u)\cap C(v)|)= \Delta + 7 + |C(u)\cap C(v)|- (d(u)- 1 + d(v)- 1)=
\Delta + 9 + |C(u)\cap C(v)|- (d(u) + d(v)) \geq\Delta + 9 +
|C(u)\cap C(v)|- 11\geq \Delta + 9 + 2 -11 =\Delta$.

Let $\{w_1, w_2, \ldots, w_{d(v)-1}\}= \{v_2, v_3, \ldots,
v_{d(v)}\}$. For any $1\leq i\leq |C(u)\cap C(v)|$, let $c(vw_i)=
i$, (i.e., $C(u)\cap C(v)= \{i | 1\leq i \leq |C(u)\cap C(v)|\}$),
$C_i= \{$$j \in C\backslash (C(u)\cup C(v))|$ $H$ contains an $(i,
j)_{(u_i, v)}$-path$\}$ and $T_i= (C\backslash (C(u)\cup
C(v)))\backslash C(w_i)$. In particular, $c(vw_1)= 1$, $c(vw_2)= 2$
and $T_1= (C\backslash (C(u)\cup C(v)))\backslash C(w_1)= \{\beta_i|
i= 1, 2, \ldots\}$, $T_2= (C\backslash (C(u)\cup C(v)))\backslash
C(w_2)= \{\alpha_i| i= 1, 2, \ldots\}$.

To complete the proof, we introduce a few symbols as defined in
\cite{BCCHM11}. A multiset is a generalized set where a member can
appear multiple times. If an element $x$ appears $t$ times in the
multiset $S$, then we say that the multiplicity of $x$ in $S$ is
$t$, and  write mult$_S(x)=t$. The cardinality of a finite multiset
$S$, denoted by $\|S\|$, is defined as $\|S\|= \sum_{x\in
S}$mult$S(x)= t$. Let $S_1$ and $S_2$ be two multisets. The join of
$S_1$ and $S_2$, denoted  $S_1\biguplus S_2$, is a multiset that
have all the members of $S_1$ as well as $S_2$. For $x\in S_1
\biguplus S_2$, mult$_{S_1 \biguplus S_2}(x) =$ mult$_{S_1}(x) +
$mult$_{S_2}(x)$. Clearly, $\|S_1 \biguplus S_2\| $ $= \|S_1\| +
\|S_2\|$.

Moreover, we define
$$S_v= \biguplus\limits_{1\leq i\leq
d(v)-1}(C(w_i)\backslash \{c(vw_i)\}).$$
 Clearly,
 $\|S_v\|=
\sum_{1\leq i\leq d(v)-1}d(w_i)- (d(v)-1)= \sum_{2\leq i\leq
d(v)}d(v_i)- (d(v)-1)$ and for any $x\in C\backslash (C(u)\cup
C(v))$, mult$_{S_v}(x)\geq 1$ by ($*_{4.1}$). Further, for any $i\in
C(u)\cap C(v)$, let
$$T'_i= \{ x\in T_i |  {\rm mult}_{S_v}(x)= 2 \},$$
$$T_0= \{x\in \biguplus\limits_{i\in C(u)\cap C(v)}T_i |
{\rm mult}_{S_v}(x)\geq 3\}.$$

\begin{claim}\label{2+}
Let $|C(u)\cap C(v)|= k + 1\geq 2$. If we can obtain an acyclic edge
coloring $c$ of $G$ when $|C(u)\cap C(v)|= k$, then for any $x\in
C\backslash (C(u)\cup C(v))$, {\rm mult}$_{S_v}(x)\geq 2$ and
$\|S_v\|\geq 2|C\backslash (C(u)\cup C(v))|$.

\end{claim}

\proof If not, then there exists $x\in C\backslash (C(u)\cup C(v))$
such that mult$S_v(x)= 1$ and $x\in C(w_1)$. Obviously, $x\in C_1$
by ($*_{4.1}$) and the $(1, x)_{(u, v)}$-path cannot pass through
$w_2$. We recolor $vw_2$ with $x$. Since mult$S_v(x)= 1$, $x\not\in
C(w_i)$, $i\neq 1$ and no bichromatic cycles are produced. Then, we
reduce the proof to the Case $|C(u)\cap C(v)|= k$. \qed

\medskip
\noindent {\bf Case 2}\ $|C(u)\cap C(v)|=2$ and $C(u)\cap C(v)= \{1,
2\}$.
\medskip

By ($*_{4.1}$), for any $j\in C\backslash (C(u)\cup C(v))$, $H$
contains an $(1, j)_{(u_1, v)}$-path or $(2, j)_{(u_2, v)}$-path and
$T_1\subseteq C_2\subseteq C(w_2)$, $T_2\subseteq C_1\subseteq
C(w_1)$, $(C\backslash (C(u)\cup C(v)))\subseteq C_1\cup C_2$. And
for any $i\in T_1$, $(2, i)_{(u_2, w_2)}$-path cannot path through
$w_1$, while for any $j\in T_2$, $(1, j)_{(u_1, w_1)}$-path cannot
pass through $w_2$. Further, by Claim \ref{2+}, $T_1\subseteq
\bigcup\limits_{3\leq i \leq d(v)-1}C(w_i)$ and $T_2\subseteq
\bigcup\limits_{3\leq i \leq d(v)-1}C(w_i)$.

\begin{claim}\label{-1}
$\|S_v\cap C(u)\| \geq |C(u)|- 1= d(u)- 2$.

\end{claim}

\proof\ If not, then there are at least two colors  $x_1, x_2\in
C(u)\backslash S_v$. First, remove the colors of $vw_1$ and $vw_2$,
then color $vw_1$ with $x_1$, $vw_2$ with $x_2$, i.e., $c(vw_1)=
x_1$, $c(vw_2)= x_2$. Since $x_1, x_2\not\in S_v$, no bichromatic
cycles are produced. Without loss of generality, we assume that
$x_1= 1$ and $x_2= 2$. First, we switch the colors of $vw_1$ and
$vw_2$. Next, if $C_1\cap C_2\neq \emptyset$, then we color $uv$
with a color in $C_1\cap C_2$. Otherwise, $C_1\cap C_2= \emptyset$.
Since $d(u_1)\leq \Delta$, $|C_1| + |C_2|= |C_1 \cup C_2|\geq
|C\backslash (C(u)\cup C(v))|\geq \Delta$ and $1\in C(u_1)$, there
exits $b_1\in C_2\backslash C(u_1)$. Then we color $uv$ with $b_1$
and no bichromatic cycles are produced by ($*_{4.1}$) and Lemma
\ref{exchange}.  \qed

\begin{claim}\label{cap}
For any $i\in C\{1, 2\}$, $(C(v)\backslash \{1, 2\})\cap C(w_i)\neq
\emptyset$.

\end{claim}

\proof If not, $(C(v)\backslash \{1, 2\})\cap C(w_{i_0})= \emptyset$
for some $i_0\in \{1, 2\}$, then we recolor $vw_{i_0}$ with a color
in $T_{i_0}$ and reduce the proof to the Case 1. \qed

Recall that $|C\backslash (C(u)\cup C(v))|\geq \Delta$. It follows
that $|T_i|\geq 2$ for any $i\in \{1, 2\}$.

\medskip
\noindent {\bf Case 2.1}\ ($A_{4.1}$) holds and $d(v)= 4$.
\medskip

Assume that $c(vw_3)= 7$. By Claim \ref{cap}, $7\in C(w_1)$ and
$7\in C(w_2)$. By Claim \ref{2+} and Claim \ref{-1}, for any $x\in
C\backslash (C(u)\cup C(v))$, mult$_{S_v}(x)\geq 2$ and $\|S_v\cap
C(u)\| \geq d(u)- 2$. And recall that $7\in C(w_1)\cap C(w_2)$, we
have $\|S_v\|\geq 2|C\backslash (C(u)\cup C(v))| + (|C(u)|- 1) +
|\{7\}| + |\{7\}|= 2(\Delta+ 11- (d(u)+ d(v))) + (d(u)- 2)+ 2=
2\Delta + 14- d(u)$. However, $\|S_v\|= \sum_{2\leq i\leq
d(v)}d(v_i)- (d(v)-1)= \sum_{2\leq i\leq d(v)}d(v_i)- 3 \leq 2\Delta
+ d(v_2)- 3$. Clearly, $2\Delta + 14- d(u) \leq 2\Delta + d(v_2)-
3$, which implies that $d(u)+ d(v_2)\geq 17$. Hence, $d(u) + d(v_2)
= 17$ and $\|S_v\cap C(u)\| = d(u)- 2$. Assume that there exists
$k\in C(u)\backslash S_v$. By ($*_{4.1}$), for any $i\in T_1$, $(2,
i)_{(u, v)}$-path cannot pass through $w_1$ since $i\not\in C(w_1)$.
If $k= 1$, recolor $vw_3$ with $1$, $vw_1$ with $\beta_1$ and color
$uv$ with $\alpha_1$. Similarly, $k \neq 2$. Otherwise, $k\in
C(u)\backslash\{1, 2\}$. Without loss of generality, assume that $k=
3$. If $H$ contains no $(3, \beta_1)_{(u, w_2)}$-path, we recolor
$vw_2$ with $3$ and color $uv$ with $\beta_1$. Otherwise, recolor
$vw_3$ with $3$, $vw_2$ with $\alpha_1$ and color $uv$ with
$\beta_1$. By ($*_{4.1}$) and Lemma \ref{exchange}, no bichromatic
cycles are produced.

\medskip
\noindent {\bf Case 2.2}\ ($A_{4.2}$) holds and $d(v)= 5$.
\medskip

Assume that $c(vw_3)= 6$ and $c(vw_4)= 7$ . By Claim \ref{cap},
there exists $a, b\in \{6, 7\}$ such that $a\in C(w_1)$ and $b\in
C(w_2)$.

By Claim \ref{2+} and Claim \ref{-1}, for any $x\in C\backslash
(C(u)\cup C(v))$, mult$_{S_v}(x)\geq 2$ and $\|S_v\cap C(u)\| \geq
d(u)- 2$.  And recall that $a\in C(w_1)\cap \{6, 7\}$, $b\in
C(w_2)\cap \{6, 7\}$, we have $\|S_v\|\geq 2|C\backslash (C(u)\cup
C(v))| + |T_0| + (d(u)- 2) + |\{a\}| + |\{b\}|= 2(\Delta+ 11- (d(u)+
d(v)))+ |T_0| + (d(u)- 2)+ 2= 2\Delta + 12 + |T_0|- d(u)$. However,
$\|S_v\|= \sum_{2\leq i\leq d(v)}d(v_i)- (d(v)-1)= \sum_{2\leq i\leq
d(v)}d(v_i)- 4 \leq 2\Delta + d(v_2) + d(v_3)- 4$. Clearly, $2\Delta
+ 12 + |T_0|- d(u) \leq 2\Delta + d(v_2) + d(v_3)- 4$, which implies
that $d(u)+ d(v_2) + d(v_3)\geq 16 + |T_0|$. Hence, $|T_0|\leq 3$
and $|T_0|= 3$ if and only if $d(u)+ d(v_2) + d(v_3)= 19$ and
$|C(u)\cap S_v|= d(u)- 2$. Clearly, there exists $j_0\in \{\alpha_1,
\alpha_2, \beta_1, \beta_2\}$ such that mult$_{S_v}(j_0)= 2$ and
assume that $j_0= \alpha_1\in T'_2$. Since $T_2 \subseteq C(w_1)$,
we can assume that $\alpha_1\in C(w_3)\backslash C(w_4)$. If $H$
contains no $(6, \alpha_1)_{(w_2, w_3)}$-path, then recolor $vw_2$
with $\alpha_1$ and reduce the proof to the Case 1. By ($*_{4.1}$)
and Lemma \ref{exchange}, no bichromatic cycles are produced.
Otherwise, we assume that:

\noindent($*_{4.2}$)\ $6\in C(w_2)$ and $H$ contains a $(6,
\alpha_1)_{(w_2, w_3)}$-path, which cannot pass through $w_4$.

\begin{claim}\label{same side}
\ \ \ {}

 {\rm (a)}\ If $|T'_2|\geq 2$, then $T'_2\subseteq
C(w_3)\backslash C(w_4)$. And for any $\alpha_i\in T'_2$, $H$
contains a $(6, \alpha_i)_{(w_2, w_3)}$-path, which cannot pass
through $w_4$.

{\rm (b)}\ If $|T'_1|\geq 1$, then there exists $j\in \{3, 4\}$ and
$k\in \{3, 4\}\backslash \{j\}$ such that $T'_1\subseteq
C(w_j)\backslash C(w_k)$. And for any $\beta_i\in T'_1$, $H$
contains a $(c(vw_j), \beta_i)_{(w_1, w_j)}$-path, which cannot pass
through $w_k$, and $c(vw_j)\in C(w_1)$.

\end{claim}

\proof To prove (a), we assume that there exists $\alpha_2\in
T'_2\backslash \{\alpha_1\}$ such that $\alpha_2\in C(w_4)\backslash
C(w_3)$. Then, we recolor $vw_4$ with $\alpha_1$, $vw_2$ with
$\alpha_2$. Since $\alpha_1\not\in C(w_2)$, $\alpha_2\not\in C(w_3)$
and by ($*_{4.2}$) and Lemma \ref{exchange}, no bichromatic cycles
are produced and reduce the proof to the Case 1. Thus,
$T'_2\subseteq C(w_3)\backslash C(w_4)$. If $H$ contains no $(6,
\alpha_i)_{(w_2, w_3)}$-path for some $\alpha_i\in T'_2$, then
recolor $vw_2$ with $\alpha_i$ and reduce the proof to the Case 1.
Similarly, we can obtain (b). \qed

By Claim \ref{same side}, $T'_2\subseteq C(w_3)$. Since $\alpha_i\in
C(w_3)\cap C(w_4)$ if mult$_{S_v}(\alpha_i)\geq 3$ for any
$\alpha_i\in T_2\backslash T'_2$, we have $T_2\subseteq C(w_3)$.
Thus, $d(w_3)\geq |T_2| + |\{6\}|= |T_2| + 1$. Similarly, there
exists $j\in \{3, 4\}$ such that $T_1\subseteq C(w_j)$ and
$d(w_j)\geq |T_1| + |\{c(vw_j)\}|= |T_1| + 1$.

\begin{claim}\label{C(u)}

$C(u)\subseteq S_v$.

\end{claim}

\proof If $2\not\in S_v$, we recolor $vw_3$ with $2$, $vw_2$ with
$\alpha_1$ and color $uv$ with $\beta_1$. By ($*_{4.1}$) and Lemma
\ref{exchange}, no bichromatic cycles are produced. Assume that
$(C(u)\backslash \{1, 2\})\backslash S_v\neq \emptyset$ and $3\in
C(u)\backslash S_v$. If $H$ contains no $(3, \beta_1)_{(u_3,
w_2)}$-path, then recolor $vw_2$ with $3$ and color $uv$ with
$\beta_1$. Otherwise, recolor $vw_3$ with $3$, $vw_2$ with
$\alpha_1$ and color $uv$ with $\beta_1$.

Assume that $1\not\in S_v$. If $|T'_1|\geq 1$, $\beta_1\in T'_1$ and
assume that $\beta_1\in C(w_3)\backslash C(w_4)$ (it is similar if
$\beta_1\in C(w_4)\backslash C(w_3)$). Then, we recolor $vw_4$,
$vw_3$, $vw_1$ with $\beta_1$, $1$ and $\beta_2$, and color $uv$
with $\alpha_1$. Clearly, no bichromatic cycles are produced by
($*_{4.1}$) and Lemma \ref{exchange}. Otherwise, $|T'_1|= 0$ and
mult$_{S_v}(\beta_i)\geq 3$ for any $\beta_i\in T_1$, which implies
that $|T_0|\geq |T_1|$. If $H$ contains no $(6, \beta_1)_{(w_1,
w_3)}$-path, then recolor $vw_4$ with $1$, $vw_1$ with $\beta_1$ and
color $uv$ with $\alpha_1$. Otherwise, $H$ contains a $(6,
\beta_1)_{(w_1, w_3)}$-path and $6\in C(w_1)$. If $H$ contains no
$(7, \beta_1)_{(w_1, w_4)}$-path, then recolor $vw_3$ with $1$,
$vw_1$ with $\beta_1$ and color $uv$ with $\alpha_1$. Otherwise, $H$
contains a $(7, \beta_1)_{(w_1, w_4)}$-path and $7\in C(w_1)$. It
follows that $1, 6, 7\in C(w_1)$ and $|T_1|\geq \Delta + 11 - (d(u)
+ d(v))- (d(w_1)- 3)= \Delta + 11 - (d(u) + d(v))- d(w_1) + 3\geq
\Delta - d(w_1) + 3\geq 3$ since $d(u)+ d(v)\leq 11$ and $d(w_1)\leq
\Delta$. Recall that $|T_1|\leq |T_0|\leq 3$. We have $|T_0|=
|T_1|=3$, $d(u)+ d(v)= 11$, $d(w_1)= \Delta$ and $d(u)+ d(v_2) +
d(v_3)= 19$, which implies that $d(u)= d(v_2)= 6$, $d(v_i)= 7$ for
$i= 3, 4, 5$ and $\Delta= 7$. And recall that $6, 7\in C(w_1)$,
$6\in C(w_2)$, we have $\|S_v\|\geq 2|C\backslash (C(u)\cup C(v))| +
|T_0| + (|C(u)|- 1) + $mult$_{S_v}(6) + $mult$_{S_v}(7)$ $=
2(\Delta+ 11- (d(u)+ d(v)))+ |T_0| + (d(u)- 2)+ 3= 2\Delta + 3 + 4 +
3= 2\times 7+ 10= 24$. However, $\|S_v\|= \sum_{2\leq i\leq
d(v)}d(v_i)- (d(v)-1) = 6 + 3\times 7- 4= 23$, a contradiction. \qed

\begin{claim}\label{676}

$G$ contains a $((1, 7)_{(u_1, w_1)})$-path or a $(2, 7)_{(u_2,
w_2)}$-path; {\rm mult}$_{S_v}(6) + $ {\rm mult}$_{S_v}(7)\geq 3$,
$7\in C(w_2)$ or $7\in C(w_1)\cap C(w_3)$ and $H$ contains a $(6,
7)_{(w_2, w_3)}$-path.

\end{claim}

\proof If $H$ contains no $(1, 7)_{(u_1, w_1)}$-path and $(2,
7)_{(u_2, w_2)}$-path, then recolor $vw_4$ with $\alpha_1$ and color
$uv$ with $7$. By ($*_{4.2}$) and Lemma \ref{exchange}, no
bichromatic cycles are produced. Otherwise, $H$ contains a $(1,
7)_{(u_1, w_1)}$-path or $(2, 7)_{(u_2, w_2)}$-path and $7\in
C(w_1)\cup C(w_2)$. If $7\in C(w_2)$, then mult$_{S_v}(6) +
$mult$_{S_v}(7)\geq \|\{6, 7, a\}\| =3$. Otherwise, $7\in
C(w_1)\backslash C(w_2)$ and $H$ contains a $(1, 7)_{(u_1,
w_1)}$-path. If $H$ contains no $(6, 7)_{(w_2, w_3)}$-path, then
recolor $vw_2$ with $7$, $vw_4$ with $\alpha_1$ and color $uv$ with
$\beta_1$. And no bichromatic cycles are produced by ($*_{4.2}$) and
Lemma \ref{exchange}. Otherwise, $H$ contains a $(6, 7)_{(w_2,
w_3)}$-path and $7\in C(w_3)$. Hence, mult$_{S_v}(6) +
$mult$_{S_v}(7)\geq$ mult$_{C(w_1)}(7) +$ mult$_{C(w_2)}(6)+
$mult$_{C(w_3)}(7)= 3$. \qed

Recall that $C(u)\subseteq S_v$ by Claim \ref{C(u)}. Hence,
$\|S_v\|\geq 2|C\backslash (C(u)\cup C(v))| + |T_0| + |C(u)| + $
mult$_{S_v}(6) + $ mult$_{S_v}(7)$ $= 2(\Delta+ 11- (d(u)+ d(v)))+
|T_0| + d(u)- 1 + $mult$_{S_v}(6) + $mult$_{S_v}(7)\geq 2\Delta+ 12-
d(u)+ |T_0| - 1 + 3= 2\Delta+ 14 + |T_0|- d(u)$. However, $\|S_v\|=
\sum_{2\leq i\leq d(v)}d(v_i)- (d(v)-1)= \sum_{2\leq i\leq
d(v)}d(v_i)- 4\leq 2\Delta + d(v_2) + d(v_3)- 4$. Hence, $2\Delta +
14 + |T_0|- d(u) \leq 2\Delta + d(v_2) + d(v_3)- 4$, which implies
that $d(u)+ d(v_2) + d(v_3)\geq 18 + |T_0|$.

Since $d(u)+ d(v_2) + d(v_3)\leq 19$, we have:

\noindent($*_{4.3}$)\ (i)\ $d(u)+ d(v_2) + d(v_3)= 18$, $|T_0|= 0$
and mult$_{S_v}(6) + $mult$_{S_v}(7)= 3$; or (ii)\ $d(u)+ d(v_2) +
d(v_3)= 19$, $|T_0|\leq 1$ and $|T_0|= 1$ if and only if
mult$_{S_v}(6) + $mult$_{S_v}(7)= 3$.

Let $\kappa_1= |C(w_1)\cap (C(u)\cup C(v))|$, $\kappa_2= |C(w_2)\cap
(C(u)\cup C(v))|$. We have $|T_1|= |(C\backslash (C(u)\cup
C(v)))\backslash C(w_1)|= \Delta + 11 -(d(u) + d(v))- (d(w_1)-
\kappa_1)\geq \Delta - d(w_1) + \kappa_1\geq \kappa_1$ and $|T_2|=
|(C\backslash (C(u)\cup C(v)))\backslash C(w_2)|= \Delta + 11 -(d(u)
+ d(v))- (d(w_2)- \kappa_2)\geq \Delta - d(w_2) + \kappa_2\geq
\kappa_2$.

\begin{claim}\label{6767}
If {\rm mult}$_{S_v}(6) + $ {\rm mult}$_{S_v}(7)\geq 4$, then

{\rm (i)}\ $|T_0|= 0$, and $d(u)= d(v_2)= 6$, $d(v_i)= 7$, $i= 3, 4,
5$;

{\rm (ii)}\ {\rm mult}$_{S_v}(6) + $ {\rm mult}$_{S_v}(7)= 4$, and
for any $i\in C(u)$, {\rm mult}$_{S_v}(i)= 1$;

{\rm (iii)}\ $|C\backslash (C(u)\cup C(v))| = \Delta$ and $|T_1|=
\Delta - d(w_1) + \kappa_1$, $|T_2|= \Delta- d(w_2) + \kappa_2$.

{\rm (iv)}\ If $|T_1|= \kappa_1$, then $\Delta= d(w_1)= 7$; If
$|T_2|= \kappa_2$, then $\Delta= d(w_2)= 7$.

\end{claim}

\proof Since $\|S_v\|\geq 2|C\backslash (C(u)\cup C(v))| + |T_0| +
$$\sum_{x\in C(u)}$mult$S_v(x)$$+ $ mult$S_v(6) + $ mult$S_v(7)$ $\geq
2|C\backslash (C(u)\cup C(v))| + |T_0| + |C(u)| + 4= 2(\Delta+ 11-
(d(u)+ d(v))) + |T_0|+ d(u)- 1 + 4\geq 2\Delta+ 12 + |T_0| -d(u)+ 3=
2\Delta + 15+ |T_0|- d(u)$ and $\|S_v\|= \sum_{2\leq i\leq
d(v)}d(v_i)- (d(v)-1)= \sum_{2\leq i\leq d(v)}d(v_i)- 4\leq 2\Delta
+ d(v_2) + d(v_3)- 4$, we have $2\Delta + 15+ |T_0|- d(u)\leq
2\Delta + d(v_2) + d(v_3)- 4$, which implies that $d(u) + d(v_2) +
d(v_3) \geq 19 + |T_0|$. Since $|T_0|\geq 0$, we have $d(u) + d(v_2)
+ d(v_3) = 19$, $|T_0|= 0$. It follows that (ii)-(iv) holds. \qed

By Claim \ref{same side}, Claim \ref{6767}, we can further assume
that:

\medskip
\noindent($*_{4.4}$)\ If mult$_{S_v}(6) + $mult$_{S_v}(7)\geq 4$,
then $T'_2= T_2\subseteq C(w_3)$, $T'_1= T_1\subseteq C(w_j)$ for
some $j\in \{3, 4\}$, $c(vw_j)\in C(w_1)$ and for any $\beta_i\in
T_1$, $H$ contains a $(c(vw_j), \beta_i)_{(w_1, w_j)}$-path, which
cannot path through $w_k$, $k\in \{3, 4\}\backslash \{j\}$.
\medskip

Since $|T_0|\leq 1$ by ($*_{4.3}$) and $|T_1|\geq 2$, we have
$|T'_1|\geq 1$ and assume that $\beta_1\in T'_1$. Now we need to
consider the following.

\medskip
\noindent {\bf Case 2.2.1}\ Assume that $\beta_1\in C(w_3)$.
\medskip

By Claim \ref{same side}, $6\in C(w_1)$. By Claim \ref{676}, $7\in
C(w_2)$ or $7\in C(w_1)\cap C(w_3)$ and $H$ contains a $(6,
7)_{(w_2, w_3)}$-path.

$\bullet$ \ $7\in C(w_1)\cap C(w_2)$.

Then mult$_{S_v}(6) + $mult$_{S_v}(7)\geq 4$. By Claim \ref{6767},
$d(w_3)\leq 7$, $6, 7\not\in C(w_3)\cup C(w_4)$, and for any $i\in
C(u)$, mult$_{S_v}(i)= 1$, $T'_1= T_1\subseteq C(w_3), T'_2=
T_2\subseteq C(w_3)$. Since $|T_1|\geq\lambda_1\geq 3$,
$|T_2|\geq\lambda_2\geq 3$, $d(w_3)\leq 7$ and $d(w_3)\geq |T_1| +
|T_2| + |\{6\}|\geq 7$, we have $T_1= \{\beta_1, \beta_2,
\beta_3\}$, $T_2= \{\alpha_1, \alpha_2, \alpha_3\}$, $C(w_3)= \{6,
\alpha_1, \alpha_2, \alpha_3, \beta_1, \beta_2, \beta_3\}$. We
recolor $vw_3$ with $2$, $vw_2$ with $\alpha_1$ and color $uv$ with
$\beta_1$. By ($*_{4.1}$) and Lemma \ref{exchange}, no bichromatic
cycles are produced.

$\bullet$ \ $7\in C(w_2)\backslash C(w_1)$.

If $H$ contains no $(6, 7)_{(w_1, w_3)}$-path, we recolor $vw_4$
with $\alpha_1$, $vw_1$ with $7$ and color $uv$ with $\alpha_2$. By
($*_{4.1}$), ($*_{4.2}$) and Lemma \ref{exchange}, no bichromatic
cycles are produced. Otherwise, $H$ contains a $(6, 7)_{(w_1,
w_3)}$-path, which implies that $7\in C(w_3)$. Thus, mult$_{S_v}(6)
+ $mult$_{S_v}(7)\geq 4$.

By Claim \ref{6767}, $d(w_3)\leq 7$, $6\not\in C(w_4)$, and for any
$i\in C(u)$, mult$_{S_v}(i)= 1$, $T'_1= T_1\subseteq C(w_3), T'_2=
T_2\subseteq C(w_3)$. Since $|T_1|\geq\kappa_1\geq 2$,
$|T_2|\geq\kappa_2\geq 3$, and $d(w_3)\geq |T_1| + |T_2| + |\{6,
7\}|\geq 7$, we have $T_1= \{8, 9\}$, $T_2= \{10, 11, 12\}$,
$C(w_3)= \{6, 7, 8, 9, 10, 11, 12\}$. Further, we may assume that
$d(w_1) = d(w_2)= d(w_3)= \Delta= 7$ and $d(w_4)= 6$ by Claim
\ref{6767}, and $C(w_1)= \{1, 6, 10, 11, 12, 13, 14\}$, $C(w_2)=
\{2, 6, 7, 8, 9, 13, 14\}$, $C(w_4)= \{1, 2, 3, 4, 5, 7\}$. If $H$
contains no $(1, 13)_{(w_1, w_4)}$-path, recolor $vw_4$ with $13$,
$vw_3$ with $2$, $vw_2$ with $10$ and color $uv$ with $8$.
Otherwise, $H$ contains a $(1, 13)_{(w_1, w_4)}$-path and $(2,
13)_{(w_2, u)}$-path by ($*_{4.1}$), which cannot pass through
$w_4$. Recolor $vw_4$, $vw_1$, $vw_3$ with $13$, $8$, $1$,
respectively, and color $uv$ with $10$. By ($*_{4.1}$) and Lemma
\ref{exchange}, no bichromatic cycles are produced.

$\bullet$ \ $7\in C(w_1)$ and $7\in C(w_3)\backslash C(w_2)$.

Note that  mult$_{S_v}(6) + $mult$_{S_v}(7)\geq 4$. By Claim
\ref{6767}, $d(w_3)\leq 7$, $6\not\in C(w_4)$, and for any $i\in
C(u)$, mult$_{S_v}(i)= 1$, $T'_1= T_1\subseteq C(w_3), T'_2=
T_2\subseteq C(w_3)$. Since $|T_1|\geq 3$, $|T_2|\geq 2$, and
$d(w_3)\geq |T_1| + |T_2| + |\{6, 7\}|\geq 7$, we have $T_1= \{10,
11, 12\}$, $T_2= \{8, 9\}$, $C(w_3)= \{6, 7, 8, 9, 10, 11, 12\}$.
Further, we may assume that $d(w_1) = d(w_2)= d(w_3)= \Delta= 7$ and
$d(w_4)= 6$ by Claim \ref{6767}, and $C(w_1)= \{1, 6, 7, 8, 9, 13,
14\}$, $C(w_2)= \{2, 6, 10, 11, 12, 13, 14\}$, $C(w_4)= \{1, 2, 3,
4, 5, 7\}$. If $H$ contains no $(1,13)_{(w_1, w_4)}$-path, recolor
$vw_4$ with $13$, $vw_3$ with $2$, $vw_2$ with $8$ and color $uv$
with $10$. Otherwise, $H$ contains a $(1, 13)_{(w_1, w_4)}$-path and
$(2, 13)_{(w_2, u)}$-path by ($*_{4.1}$), which cannot pass through
$w_4$. Recolor $vw_4$, $vw_1$, $vw_3$ with $13$, $10$, $1$,
respectively, and color $uv$ with $8$. By ($*_{4.1}$) and Lemma
\ref{exchange}, no bichromatic cycles are produced.

\medskip
\noindent {\bf Case 2.2.2}\ Assume that $\beta_1\in C(w_4)$.
\medskip

By Claims \ref{same side} and \ref{676}, $7\in C(w_1)$ and if
$7\not\in C(w_2)$, then $7\in C(w_3)$, and $H$ contains a $(6,
7)_{(w_2, w_3)}$-path. If $6\not\in C(w_1)$ and $G$ contains no $(2,
6)_{(u_2, w_2)}$-path, then recolor $ww_3$ with $\beta_1$ and color
$uv$ with $6$. Otherwise, $G$ contains a $(2, 6)_{(u_2, w_2)}$-path
if $6\not\in C(w_1)$. If $6\not\in C(w_1)$ and $H$ contains no $(6,
7)_{(w_1, w_4)}$-path, then recolor $vw_1$ with $6$, $vw_3$ with
$\beta_1$ and color $uv$ with $\alpha_1$. By ($*_{4.4}$) and Lemma
\ref{exchange}, no bichromatic cycles are produced. Otherwise, we
assume that:

\medskip
\noindent($*_{4.5}$)\ (i)\ $7\in C(w_2)$, or $7\in C(w_3)$ and $H$
contains a $(6, 7)_{(w_2, w_3)}$-path. (ii)\  $6\in C(w_1)$, or
$6\in C(w_4)$ and $H$ contains a $(6, 7)_{(w_1, w_4)}$-path.

\medskip

By ($*_{4.5}$), we have $7\in C(w_2)\cup C(w_3)$, $6\in C(w_1)\cup
C(w_4)$ and mult$_{S_v}(6) + $mult$_{S_v}(7)\geq 4$. By Claim
\ref{6767} and ($*_{4.3}$), ($*_{4.4}$), mult$_{S_v}(6) +
$mult$_{S_v}(7)= 4$, $d(w_i)\leq 7$ for any $1\leq i\leq 4$, and for
any $j\in C(u)$, mult$_{S_v}(j)= 1$. By the previous analysis, we
need to consider the following subcase.

$\bullet$ \ $c(vv_5)= 1$ and $7\in C(v_5)$.

If $c(uv_5)= 2$, then we color $uv$ with a color in $\{8, 9, \ldots,
\Delta + 7\}\backslash C(v_5)$. Otherwise, $c(uv_5)= 3$. Assume that
$c(vw_2)= 2$ and $c(vw_4)= 7$.

If there is no $(1, 3)_{(v, w_2)}$-path in $H$, recolor $vw_2$ with
$3$ and color $uv$ with a color in $\{8, 9, \ldots, \Delta +
7\}\backslash C(v_5)$. Since mult$_{S_v}(3)= 1$ and $3\in C(v_5)$,
there is no bichromatic cycles. Otherwise, $H$ contains a $(1,
3)_{(v, w_2)}$-path and $1\in C(w_2)$. If for some color $i\in \{8,
9, \ldots, \Delta + 7\}\backslash C(v_5)$, there is no an $(i,
j)_{(u, v_5)}$-path for any $j\in \{4, 5\}$, then we recolor $uv_5$
with $i$ and color $uv$ with $3$. Otherwise, for each color $i\in
\{8, 9, \ldots, \Delta+ 7\}\backslash C(v_5$, $H$ contains an $(i,
j)_{(u, v_5)}$-path for some color $j\in \{4, 5\}$. Without loss of
generality, assume that $4\in C(v_5)$ and $H$ contains a $(\beta_1,
4)_{(u, v_5)}$-path. Since $d(v_5)= 7$ and $d(w_2)\leq 7$, it
suffice to assume that $C(v_5)= \{1, 3, 4, 7, 8, 9, 10\}$, $C(w_2)=
\{1, 2, 6, 11, 12, 13, 14\}$ and $\Delta = 7$. Then recolor $vv_5$,
$vw_2$, $vw_4$ with $\beta_1, 4, 1$ and color $uv$ with $8$. By the
previous assumption and Lemma \ref{exchange}, no bichromatic cycles
are produced.

$\bullet$ \ $c(vv_5)= 6$, $\beta_1\not\in C(v_5)$ and $c(vw_1)= 1$,
$c(vw_2)= 2$, $v(vw_4)= 7$

If $c(uv_5)\in \{3, 4, 5\}$, then we recolor $vw_2$ with $c(uv_5)$
and color $uv$ with $\beta_1$. Since mult$_{S_v}(3)= 1$, $3\not\in
C(v_5)$ and no bichormtic cycles are produced. Otherwise,
$c(uv_5)\in \{1, 2\}$. If for some color $i\in \{8, 9, \ldots,
\Delta + 7\}\backslash C(v_5)$, there is no an $(i, j)_{(u,
v_5)}$-path for any $j\in \{3, 4, 5\}$, then we recolor $uv_5$ with
$i$ and reduce the proof to the Case 1. Otherwise, for each color
$i\in \{8, 9, \ldots, \Delta + 7\}\backslash C(v_5)$, $H$ contains
$(i, j)_{(u, v_5)}$-path for some color $j\in \{3, 4, 5\}$. Without
loss of generality, assume that $3\in C(v_5)$ and $H$ contains a
$(\beta_1, 3)_{(u, v_5)}$-path.

Assume that $c(uv_5)= 1$. Since mult$_{S_v}(3)= $mult$_{S_v}(1)= 1$,
$1, 3\not\in (C(w_1)\backslash \{1\})\cup C(w_2)$. If $H$ contains
no $(3, 6)_{(v_5, w_1)}$-path, then we recolor $vw_1$ with $3$,
$vw_2$ with $1$ and color $uv$ with $\beta_1$. Otherwise, $H$
contains a $(3, 6)_{(v_5, w_1)}$-path. Recolor $vw_2$ with $3$,
$vv_5$ with $\beta_1$ and color $uv$ with $6$. Since $H$ contains no
$(7, \beta_1)_{(w_1, w_4)}$-path by Claim \ref{same side}, no
bichromatic cycles are produced even $7\in C(v_5)$ by Lemma
\ref{exchange}. Otherwise, assume that $c(uv_5)= 2$. Since
mult$_{S_v}(3)= $mult$_{S_v}(2)= 1$, $2, 3\not\in C(w_1)$. If $H$
contains no $(3, 6)_{(v_5, w_1)}$-path, then we recolor $vw_1$ with
$3$ and color $uv$ with $\beta_1$. Otherwise, $H$ contains a $(3,
6)_{(v_5, w_1)}$-path and $6\in C(w_1)$. Recolor $vw_1$ with $3$,
$vv_5$ with $\beta_1$ and color $uv$ with $6$. Since $H$ contains no
$(7, \beta_1)_{(w_1, w_4)}$-path by Claim \ref{same side}, no
bichromatic cycles are produced even $7\in C(v_5)$ by Lemma
\ref{exchange}.

\medskip
\noindent {\bf Case 3}\ $|C(u)\cap C(v)|= 3$ and $C(u)\cap C(v)=
\{1, 2, 3\}$.
\medskip

Assume that $c(vw_3)= 3$. By ($*_{4.1}$), for any $j\in C\backslash
(C(u)\cup C(v))$, $H$ contains an $(i, j)_{(u, v)}$-path for some
$i\in \{1, 2, 3\}$ and $T_1\subseteq C_2\cup C_3\subseteq C(w_2)\cup
C(w_3)$, $T_2\subseteq C_1\cup C_3\subseteq C(w_1)\cup C(w_3)$, and
$T_3\subseteq C_1\cup C_2\subseteq C(w_1)\cup C(w_2)$. It is easy to
derive that $|C\backslash (C(u)\cup C(v))|= \Delta + 7- |C(u)\cup
C(v)|= \Delta + 7 + 2 + 3- (d(u) + d(v)) = \Delta + 12- (d(u) +
d(v)) \geq \Delta +1$, $|T_i|= |(C\backslash (C(u)\cup C(v)))
\backslash C(w_i)|= |C\backslash (C(u)\cup C(v))|- (d(w_i)-1)\geq
2$, $i= 1, 2, 3$.

\medskip
\noindent {\bf Case 3.1}\ ($A_{4.1}$) holds and $d(v)= 4$.
\medskip

Assume that $1\not\in S_v$. If $H$ contains no $(2, \beta_1)_{(w_1,
w_2)}$-path, then we recolor $vw_3$ with $1$, $vw_1$ with $\beta_1$
and reduce the proof to the Case 2. Otherwise, $H$ contains a $(2,
\beta_1)_{(w_1, w_2)}$-path and a $(3, \beta_1)_{(u, w_3)}$-path by
($*_{4.1}$) and Lemma \ref{exchange}. Then we recolor $vw_1$ with
$\beta_1$, $vw_2$ with $1$ and reduce the proof to the Case 2.
Otherwise, $1\in S_v$ and $2, 3\in S_v$. Further, if $C(u)\backslash
S_v\neq\emptyset$ and $4\in C(u)\backslash S_v$, we recolor $vw_1$
with $4$ and then   $4\in S_v$ similarly. Now we assume that
$C(u)\subseteq S_v$. By Claim \ref{2+}, for any $x\in C\backslash
(C(u)\cup C(v))$, mult$_{S_v}(x)\geq 2$ and $\|S_v\|\geq
2|C\backslash (C(u)\cup C(v))|$. So $|S_v|\geq 2|C\backslash
(C(u)\cup C(v))| + |C(u)|= 2(\Delta + 7 -(d(u)- 1+ 3- 3))+ (d(u)-
1)= 2\Delta + 15- d(u)$. $\|S_v\|= \sum_{v_i\neq u\in N(v)}d(v_i)-
(d(v)-1)= \sum_{v_i\neq u\in N(v)}d(v_i) - 3 \leq 2\Delta + d(v_2)-
3$. Hence, $2\Delta + 15- d(u)\leq 2\Delta + d(v_2)- 3$, which
implies that $d(u)+ d(v_2)\geq 18$, a contradiction.

\medskip
\noindent {\bf Case 3.2}\ ($A_{4.2}$) holds and $d(v)= 5$.
\medskip

Assume that $c(w_4)= 6$. Recall that $|T_i|\geq 2$ for any $i\in
\{1, 2, 3\}$.

\begin{claim}\label{abc}

For $i\in \{1, 2, 3\}$, $C(w_i)\cap C(v)\neq \{i\}$.

\end{claim}

\proof If $C(w_1)\cap C(v)= \{1\}$, then we recolor $vw_1$ with
$\beta_1$ and reduce the proof to the Case 2. Hence,  $C(w_1)\cap
C(v)\neq \{1\}$ and assume that $\{1, a\}\subseteq C(v)\cap C(w_1)$.
Similarly, we assume that $\{2, b\}\subseteq C(v)\cap C(w_2)$, $\{3,
c\}\subseteq C(v)\cap C(w_3)$, where $a, b, c\in C(v)$. \qed

\begin{claim}\label{23}

For $i= 1, 2, 3$, $T_i\subseteq C(w_j)$ for any $j\in \{1, 2,
3\}\backslash \{i\}$ and $T_i\cap T_j= \emptyset$ for any $i, j\in
\{1, 2, 3\}$.

\end{claim}

\proof For any $\beta_i\in T_1\backslash T'_1$, it is clear that
$\beta_i\in C(w_2)\cap C(w_3)$. If $T'_1\neq\emptyset$ and assume
that $\beta_1\in T'_1\backslash C(w_3)$, then $H$ contains a $(2,
\beta_1)_{(u, w_2)}$-path by ($*_{4.1}$) and $\beta_1\in C(w_4)$
since mult$_{S_v}(\beta_1)\geq 2$ by Claim \ref{2+}. If $H$ contains
no $(\beta_1, 6)_{(w_1, w_4)}$-path, we recolor $vw_1$ with
$\beta_1$. Otherwise, $H$ contains a $(\beta_1, 6)_{(w_1,
w_4)}$-path, recolor $vw_3$ with $\beta_1$. Clearly, no bichromatic
cycles are produced by Lemma \ref{exchange} and reduce the proof to
the Case 2. Hence, $T'_1\subseteq T_1\subseteq C(w_2)\cap C(w_3)$
and it is similar for $i= 2, 3$. \qed

By Claims \ref{abc} and  \ref{23}, $|T_i|\geq 3$ and $d(w_1)\geq
|\{1, a\}| + |T_2|+ |T_3|$, $d(w_2)\geq |\{2, b\}| + |T_1|+ |T_3|$,
$d(w_3)\geq |\{3, c\}| + |T_1|+ |T_2|$.

Let $d(u)= 6- k$, $0\leq k\leq 2$. Then $|C\backslash (C(u)\cup
C(v))|= \Delta + 7 - (d(u)- 1+ d(v)- 1) + |C(u)\cap C(v)|= \Delta +
1 + k$, $|T_i|\geq \Delta + 1 + k- (\Delta - 2)\geq k + 3$ and
$d(w_i)\geq 2+ 2(k + 3) = 2k + 8$ for $i\in \{1, 2, 3\}$. Noting
that $d(v_2)\geq d(u)$, we have $d(v_2) + d(v_3) \geq 2k +8 +
\min\{6- k, 2k + 8\}= 2k + 8+ 6- k = k + 14$ and $d(u)+ d(v_2) +
d(v_3) \geq 6- k + k +14 = 20$, a contradiction. \qed

\medskip
\noindent {\bf Case 4}\ $|C(u)\cap C(v)|= 4$ and $C(u)\cap C(v)=
\{1, 2, 3, 4\}$.
\medskip

Then $d(v)= 5$ and assume that $c(vw_3)= 3$, $c(vw_4)= 4$. Clearly,
$|C\backslash (C(u)\cup C(v))|\geq \Delta +2$ and for any $i\in \{1,
2, 3, 4\}$, $|T_i|\geq 3$.

If $C(w_i)\cap C(v)= \{i\}$ for some $i\in \{1, 2, 3, 4\}$, we
recolor $vw_i$ with a color in $T_i$ and reduce the proof to the
Case 3. Otherwise, for any $i\in \{1, 2, 3, 4\}$, $C(w_i)\cap
C(v)\neq \{i\}$. Assume that $\{1, a\}\subseteq C(w_1)\cap C(v)$,
$\{2, b\}\subseteq C(w_2)\cap C(v)$, $\{3, c\}\subseteq C(w_3)\cap
C(v)$, $\{4, d\}\subseteq C(w_4)\cap C(v)$.

If $d(u)= 5$, then $|C\backslash (C(u)\cup C(v))|\geq \Delta + 3$
and $d(v_2)+ d(v_3)\leq 18- 5= 13$. Hence, $2\Delta+ 6 + 4\leq
2|C\backslash (C(u)\cup C(v))| + |\{a, b, c, d\}|\leq\|S_v\|=
\sum_{2\leq i\leq d(v)}d(v_i)- (d(v)-1)= \sum_{2\leq i\leq
d(v)}d(v_i)- 4 \leq 2\Delta + d(v_2)+ d(v_3)- 4\leq 2\Delta + 13- 4=
2\Delta + 9$, a contradiction.

If $d(u)= 6$, then $|C\backslash (C(u)\cup C(v))|\geq \Delta + 2$
and $d(v_2)+ d(v_3)\leq 19- 6= 13$. Hence, $2\Delta+ 4 + 4 +
|T_0|\leq 2|C\backslash (C(u)\cup C(v))| + |\{a, b, c, d\}|+
|T_0|\leq \|S_v\|= \sum_{2\leq i\leq d(v)}d(v_i)- (d(v)-1)=
\sum_{2\leq i\leq d(v)}d(v_i)- 4 \leq 2\Delta + d(v_2)+ d(v_3)-
4\leq 2\Delta + 13- 4= 2\Delta + 9$, i.e., $2\Delta+ 8 + |T_0|\leq
2\Delta+ 9$, which implies that $|T_0|\leq 1$ and $|T'_1|\geq 1$
since $|T_1|\geq 4$. Assume that $\beta_1\in T'_1$ and $\beta_1\in
(C(w_2)\cap C(w_3))\backslash C(w_4)$. By ($*_{4.1}$), $H$ contains
a $(i, \beta_1)_{(u, v)}$-path for some $i\in \{2, 3\}$. Without
loss of generality, assume that $H$ contains a $(2, \beta_1)_{(u,
v)}$-path. If $H$ contains no $(3, \beta_1)_{(w_1, w_3)}$-path, we
recolor $vw_1$ with $\beta_1$ and reduce the proof to the Case 3.
Otherwise, recolor $vw_4$ with $\beta_1$ and reduce the proof to the
Case 3.\qed

%%%%%%%%%%%%%%%%%%%%%%%%%%%%%%%%%%%%%%%%%%%%%%%%%%%%%%%%%%%%%%%%%%%%%%%%%%%

\end{document}